\documentclass[12pt]{amsart}
\usepackage{amsmath}
\usepackage{amssymb}
\usepackage{latexsym}

\newcommand{\Zint}{\mathbb {Z}}    
     
\newcommand{\Cplx}{\mathbb {C}}     

\newcommand{\halmos}{\rule{5pt}{5pt}}

\numberwithin{equation}{section}

\newtheorem{prop}{\bf Proposition}[section]

\newtheorem{thm}[prop]{\bf Theorem}
\newtheorem{lemma}[prop]{\bf Lemma}

\newtheorem{exa}{\bf Example}

\begin{document}

\title[Hermite-Krichever Ansatz]
{Fuchsian equation, Hermite-Krichever Ansatz and Painlev\'e equation}
\author{Kouichi Takemura}
\address{Department of Mathematical Sciences, Yokohama City University, 22-2 Seto, Kanazawa-ku, Yokohama 236-0027, Japan.}
\email{takemura@yokohama-cu.ac.jp}


\begin{abstract}
Several results on Heun's equation are generalized to a certain class of Fuchsian differential equations.
Namely, we obtain integral representations of solutions and develop Hermite-Krichever Ansatz on them.
In particular, we investigate linear differential equations that produce Painlev\'e equation by monodromy preserving deformation and obtain solutions of the sixth Painlev\'e equation which include Hitchin's solution.
The relationship with finite-gap potential is also discussed.
\end{abstract}

\maketitle

\section{Introduction}

It is well known that a Fuchsian differential equation with three singularities is transformed to a Gauss hypergeometric equation, and plays important roles in substantial fields in mathematics and physics. Several properties of solutions to the hypergeometric equation have been explained in various textbooks.

A canonical form of a Fuchsian equation with four singularities is written as 
\begin{equation}
\left( \! \left(\frac{d}{dw}\right) ^2 \! + \left( \frac{\gamma}{w}+\frac{\delta }{w-1}+\frac{\epsilon}{w-t}\right) \frac{d}{dw} +\frac{\alpha \beta w -q}{w(w-1)(w-t)} \right)\tilde{f}(w)=0
\label{Heun}
\end{equation}
with the condition 
\begin{equation}
\gamma +\delta +\epsilon =\alpha +\beta +1,
\label{Heuncond}
\end{equation}
and is called Heun's equation.
Despite that Heun's equation was resolved in the 19th century; several results of solutions have only been recently revealed.
Namely, integral representations of solutions, global monodromy in terms of hyperelliptic integrals, relationships with the theory of finite-gap potential and the Hermite-Krichever Ansatz for the case $\gamma, \delta, \epsilon , \alpha -\beta \in \Zint +\frac{1}{2}$ are contemporary (see \cite{BE,GW,Smi,Tak1,Tak2,Tak3,Tak4,TV} etc.), though they are not written in a textbook on Heun's equation \cite{Ron}.

In this paper, we consider differential equations which have additional apparent singularities to Heun's equation.
More precisely, we consider the equation
\begin{align}
& \left\{ \frac{d^2}{dw^2}+\left( \frac{\frac{1}{2}-l_1}{w}+  \frac{\frac{1}{2}-l_2}{w-1}+  \frac{\frac{1}{2}-l_3}{w-t}+ \sum _{{i'}=1}^M \frac{-r_{i'}}{w-\tilde{b}_{i'}} \right) \frac{d}{dw} \right. \label{Feqintro} \\
& \left.  +\frac{(\sum _{i=0}^3 l_i + \sum _{{i'}=1}^M r_{i'})(-1-l_0 +\sum _{i=1}^3 l_i + \sum _{{i'}=1}^M r_{i'})w+\tilde{p}+ \sum_{{i'}=1}^M \frac{\tilde{o}_{i'}}{w-\tilde{b}_{i'}}}{4w(w-1)(w-t)}\right\}\tilde{f}(w) =0 , \nonumber
\end{align}
for the case $l_i \in \Zint _{\geq 0}$ $(0 \leq i \leq 3)$, $r_{i'} \in \Zint _{>0}$ $(1 \leq i' \leq M)$ and the regular singular points $\tilde{b}_{i'}$ $(1 \leq i' \leq M)$ are apparent.

By a certain transformation, Eq.(\ref{Feqintro}) is rewritten in terms of elliptic functions such as
\begin{align}
 & \left\{ -\frac{d^2}{dx^2} + \sum_{i=0}^3 l_i(l_i+1)\wp (x+\omega_i) \right. \label{ellDE} \\
 & \left. + \sum_{i'=1}^M \left( \frac{r_{i'}}{2} \left( \frac{r_{i'}}{2}+1\right) (\wp (x-\delta _{i'}) + \wp (x+\delta _{i'})) +\frac{s_{i'}}{\wp (x) -\wp (\delta _{i'})} \right) -E \right\} f(x) =0, \nonumber
\end{align}
with the condition that logarithmic solutions around the singularities $x= \pm \delta _{i'}$ $(i'=1,\dots ,M)$ disappear.
We then establish that solutions to Eq.(\ref{ellDE}) have an integral representation and they are also written as a form of the Hermite-Krichever Ansatz. For details see Proposition \ref{prop:Linteg} and Theorem \ref{thm:alpha}.
Note that the results on the Hermite-Krichever Ansatz are related to Picard's theorem on differential equations with coefficients of elliptic functions \cite[\S 15.6]{Inc}.
By the Hermite-Krichever Ansatz, we can obtain information on the monodromy of solutions to differential equations.

Results on integral representation and the Hermite-Krichever Ansatz are applied for particular cases.
One example is Painlev\'e equation.
For the case $M=1$ and $r_1=1$, it is known that Eq.(\ref{Feqintro}) produces the sixth Painlev\'e equation by monodromy preserving deformation (see \cite{IKSY}).
On the other hand, solutions to Eq.(\ref{ellDE}) are expressed as  a form of the Hermite-Krichever Ansatz for the case $l_i \in \Zint _{\geq 0}$ $(i=0,1,2,3)$, and we obtain an expression of monodromy.
Fixing monodromy corresponds to the monodromy preserving deformation; thus, we obtain solutions to the sixth Painlev\'e equation by fixing monodromy (see section \ref{sec:P6}). For the case $l_0=l_1=l_2=l_3=0$, we recover Hitchin's solution \cite{Hit}.

Results on integral representation and the Hermite-Krichever Ansatz are also applicable to differential equations related with finite-gap potential.
The potential of the Schr\"odinger operator as Eq.(\ref{ellDE}) for the case $M=0$ is called Treibich-Verdier potential \cite{TV}, and is an example of a finite-gap potential. For this case, the differential equation is transformed to Heun's equation. 
If $M=1$, $r_1=2$, $s_1=0$ and $b_1$ satisfies a certain algebraic equation ($s_1$ and $b_1$ appear in Eq.(\ref{ellDE})), then it is seen \cite{Tre,Smi2} that the potential is a finite-gap and is also Picard's in the sense of \cite{GW2}. 
For this potential, in this paper we provide a viewpoint from a Fuchsian equation with an apparent singularity, and more results are produced in \cite{TakF}.

This paper is organized as follows. In section \ref{sec:FDE}, we introduce Fuchsian differential equations and rewrite them to the form of elliptic functions. 
The definition of apparent singularity and its property are mentioned.
In section \ref{sec:HK}, we obtain integral representations of solutions to the differential equation of the class mentioned above and rewrite them to the form of the Hermite-Krichever Ansatz. To obtain an integral representation, we introduce doubly-periodic functions that satisfy a differential equation of order three. Some properties related with this doubly-periodic function are investigated, and we obtain another expression of solutions that looks like the form of the Bethe Ansatz (see Proposition \ref{prop:BA}).
In section \ref{sec:P6}, we consider the relationship with the sixth Painlev\'e equation. We show that solutions of the sixth Painlev\'e equation are obtained from solutions expressed in the form of the Hermite-Krichever Ansatz of linear differential equations considered in section \ref{sec:HK} by fixing monodromy. Some explicit solutions that include Hitchin's solution are displayed.
In section \ref{sec:FG}, we discuss the relationship with the results on finite-gap potential.
In section \ref{sec:rmk}, we give concluding remarks and present an open problem.
In the appendix, we note definitions and formulae for elliptic functions.

\section{Fuchsian differential equation} \label{sec:FDE}

To begin with, we introduce the following differential equation;
\begin{align}
& \left\{ \frac{d^2}{dz^2}+\left( \sum _{i=1}^3 \frac{\frac{1}{2}-l_i}{z-e_i}+ \sum _{{i'}=1}^M \frac{-r_{i'}}{z-b_{i'}} \right) \frac{d}{dz} +\frac{N(N-2l_0-1)z+p+ \sum_{{i'}=1}^M \frac{o_{i'}}{z-b_{i'}}}{4(z-e_1)(z-e_2)(z-e_3)}\right\}\tilde{f}(z) =0 ,
\label{Feq} 
\end{align}
where $N=\sum _{i=0}^3 l_i + \sum _{{i'}=1}^M r_{i'}$. This equation is Fuchsian, i.e., all singularities $\{ e_i \} _{i=1,2,3 }$, $\{ b_{i'} \} _{ i' =1, \dots ,M}$ and $\infty$ are regular. The exponents at $z=e_i$ $(i=1,2,3)$ (resp. $z=b_{i'} $ $(i'=1,\dots ,M)$) are $0$ and $l_i +1/2$ (resp. $0$ and $r_{i'} +1$), and the exponents at $z= \infty$ are $N/2$ and $(N-2l_0-1)/2$. Conversely, any Fuchsian differential equation that has regular singularities at $\{ e_i \} _{i=1,2,3 }$, $\{ b_{i'} \} _{ i' =1, \dots ,M}$ and $\infty$ such that one of the exponents at $e_i$ and $b_{i'} $ for all $i\in \{1,2,3 \}$ and $i' \in \{1, \dots ,M\}$ are zero is written as Eq.(\ref{Feq}).
By the transformation $z \rightarrow z +\alpha $, we can change to the case $e_1 +e_2 +e_3 =0$. In this paper we restrict discussion to the case $e_1 +e_2 +e_3 =0$.
We remark that any Fuchsian equation with $M+4$ singularities is transformed to Eq.(\ref{Feq}) with the condition $e_1 +e_2 +e_3 =0$.

It is known that, if $e_1 +e_2 +e_3 =0$ and $e_1 \neq e_2 \neq e_3 \neq e_1$, then there exists some periods $(2\omega _1 ,2\omega _3)$ such that $\wp (\omega _1)= e_1$ and $\wp (\omega _3 )= e_3$, where $\wp (x)$ is the Weierstrass $\wp$-function with periods $(2\omega_1, 2\omega_3)$.
We set $\omega _0 =0$ and $\omega _2 =-\omega _1 -\omega _3$. Then we have $\wp (\omega _2 )= e_2$.

Now we rewrite Eq.(\ref{Feq}) in an elliptic form. We set
\begin{equation}
\Phi (z)= \prod_{i=1}^3 (z-e_i)^{-l_i/2} \prod_{{i'}=1}^M (z-b_{i'})^{-r_{i'}/2} , \quad z= \wp(x),
\end{equation}
and $\tilde{f}(z) \Phi (z)=f(x)$. Then we have
\begin{equation}
(H-E) f(x)= 0,
\label{eq:H}
\end{equation}
where $H$ is a differential operator defined by
\begin{align}
 H= & -\frac{d^2}{dx^2} + v(x), \label{Ino} \\
v(x) = & \sum_{i=0}^3 l_i(l_i+1)\wp (x+\omega_i) \\
& + \sum_{i'=1}^M \frac{r_{i'}}{2} \left( \frac{r_{i'}}{2}+1\right) (\wp (x-\delta _{i'}) + \wp (x+\delta _{i'})) +\frac{s_{i'}}{\wp (x) -\wp (\delta _{i'})}, \nonumber
\end{align}
and
\begin{align}
& \wp (\delta _{i'})= b_{i'} , \quad ({i'}=1,\dots ,M) ,\\
& o_{i'} = -s_{i'} + r_{i'} \left\{ \frac{1}{8}r _{i'} (12b_{i'} ^2-g_2)   +\frac{1}{2}(4b_{i'} ^3- g_2 b_{i'} -g_3)\left(\sum_{i'' \neq {i'}}\frac{r_{i''}}{(b_{i'} -b_{i''})}\right) \right. \\
& \quad \quad \left. +2(l_1(b_{i'} -e_2)(b_{i'} -e_3)+l_2(b_{i'} -e_1)(b_{i'} -e_3)+l_3(b_{i'} -e_1)(b_{i'} -e_2))\right\} ,\nonumber \\
& p= E+(e_1l_1^2+e_2l_2^2+e_3l_3^2)  -2(l_1l_2e_3+l_2l_3e_1+l_3l_1e_2) -\frac{1}{2}\sum _{{i'}=1}^M b_{i'} r_{i'} ^2\\
&  \quad \quad  +2\sum _{{i'}=1}^M \sum _{i=1}^3 l_ir_{i'} (e_i+b_{i'} )+2\left(\sum _{{i'}=1}^M b_{i'} r_{i'} \right)\left(\sum _{{i'}=1}^M r_{i'} \right) ,\nonumber \\
& g_2=-4(e_1e_2+e_2e_3+e_3e_1), \quad g_3=4e_1e_2e_3. 
\end{align}
Conversely, Eq.(\ref{Feq}) is obtained from Eq.(\ref{eq:H}) by the transformation above.

We consider another expression. Set
\begin{align}
&  H_g= -\frac{d^2}{dx^2} + \sum_{i'=1}^M  \frac{r_{i'} \wp ' (x)}{\wp (x) -\wp (\delta _{i'})} \frac{d}{dx} + \left(l_0 + \sum_{i'=1}^M r_{i'}\right) \left(l_0 +1-  \sum_{i'=1}^M r_{i'}\right) \wp (x) \\
& \quad \quad +\sum_{i=1}^3 l_i(l_i+1) \wp (x+\omega_i) + \sum_{i'=1}^M \frac{\tilde{s}_{i'}}{\wp (x) -\wp (\delta _{i'})} ,\nonumber \\
& f_g(x) = f(x) \Psi _g (x), \quad \Psi _g (x)=\prod _{i'=1}^M (\wp (x) -\wp (\delta _{i'}))^{r_{i'} /2}. 
\end{align}
Then Eq.(\ref{eq:H}) is also equivalent to 
\begin{equation}
(H_g -E-C_g) f_g(x) =0,
\label{eq:Hg}
\end{equation}
where 
\begin{align}
& \tilde{s}_{i'} = s_{i'} -r_{i'} \left\{ \frac{1}{8} r _{i'} (12 b_{i'} ^2-g_2) +\frac{1}{2}(4b_{i'} ^3- g_2 b_{i'} -g_3)\left(\sum_{i'' \neq {i'}}\frac{r_{i''}}{(b_{i'} -b_{i''})}\right) \right\} ,\\
& C_g =- \frac{1}{2}\sum _{{i'}=1}^M b_{i'} r_{i'} ^2 +2\left(\sum _{{i'}=1}^M b_{i'} r_{i'} \right)\left(\sum _{{i'}=1}^M r_{i'} \right) .
\end{align}

In this paper, we consider solutions to Eq.(\ref{Feq}), which is equivalent to Eq.(\ref{eq:H}) or Eq.(\ref{eq:Hg}) for the case $l_i \in \Zint $, and the regular singular point $z= b_{i'}$ is apparent for all ${i'}$. Here, a regular singular point $x=a$ of a linear differential equation of order two is said to be apparent, if and only if the differential equation does not have a logarithmic solution at $x=a$ and the exponents at $x=a$ are integers. It is known that the regular singular point $x=a$ is apparent, if and only if the monodromy matrix around $x=a$ is a unit matrix.
Note that Smirnov investigated solutions in \cite{Smi2} with the assumptions $s_{i'}=0$ and $r_{i'} \in 2\Zint $ for all ${i'}$.

We consider the condition that the regular singular point $x=a$ is apparent.
More precisely, we describe the condition that a differential equation of order two does not have logarithmic solutions at a regular singular point $x=a$ for the case $\alpha _2 -\alpha _1 \in \Zint$, where $\alpha _1$ and  $\alpha _2$ are exponents at $x=a$.
If $\alpha _1 = \alpha _2$, then the differential equation has logarithmic solutions at $x=a$. We assume that the exponents satisfy $\alpha _2 -\alpha _1 =n \in \Zint _{\geq 1}$.
Since the point $x=a$ is a regular singular, the differential equation is written as
\begin{equation}
\left\{ \frac{d^2}{dx^2} + \sum _{j=0}^{\infty} p_j (x-a)^{j-1} \frac{d}{dx} + \sum _{j=0}^{\infty} q_j (x-a)^{j-2} \right\} f(x)=0. \label{eq:Feqxa}
\end{equation}
Let $F(t)$ be the characteristic polynomial at the regular point $x=a$, i.e. $F(t)= t^2+(p_0-1)t+q_0$. From the definition of exponents, we have $F(\alpha _1 )=F(\alpha _1  +n) =0$. 
We can now calculate solutions to Eq.(\ref{eq:Feqxa}) in the form 
\begin{equation}
f(x) = \sum _{j=0}^{\infty} c_j (x-a)^{\alpha _1 + j},
\end{equation}
where $f(x)$ is normalized to satisfy $c_0=1$. By substituting it into Eq.(\ref{eq:Feqxa}) and comparing the coefficients of $(x-a)^{\alpha _1 + j-2}$, we obtain the relations 
\begin{equation}
F(\alpha _1 +j) c_j + \sum_{j'=0}^{j-1} \{ (\alpha _1 +j') p_{j-j'} +q_{j-j'}\} c_{j'} =0. \label{eq:recj}
\end{equation}
If the positive integer $j$ satisfies $F(\alpha _1 +j) \neq 0$ (i.e. $j\neq 0,n$), then the coefficient $c_j$ is determined recursively.
For the case $j=n$, we have $F(\alpha _1 +n)=0$ and 
\begin{equation}
\sum_{j'=0}^{n-1} \{ (\alpha _1 +j') p_{n-j'} +q_{n-j'}\} c_{j'} =0. \label{eq:recn}
\end{equation}
Eq.(\ref{eq:recn}) with recursive relations (\ref{eq:recj}) for $j=1, \dots ,n-1$ is a necessary and sufficient condition that Eq.(\ref{eq:Feqxa}) does not have a logarithmic solution for the case $\alpha _2 -\alpha _1 =n \in \Zint _{\geq 1}$. In fact, if $p_0 ,q_0 , \dots , p_n ,q_n$ satisfy Eq.(\ref{eq:recn}), then there exist solutions to Eq.(\ref{eq:Feqxa}) that include two parameters $c_0$ and $c_n$. Thus any solutions are not logarithmic at $x=a$. Conversely, if Eq.(\ref{eq:recn}) is not satisfied, there exists a logarithmic solution written as $f(x)= \sum _{j=0}^{\infty} c_j (x-a)^{\alpha _1 + j} + \log (x-a) \sum _{j=n}^{\infty} \tilde{c}_j (x-a)^{\alpha _1 + j}$.

It follows from $\wp (\delta _{i'})= b_{i'}$ and $\wp '(\delta _{i'}) \neq 0$ that, the monodromy matrix to Eq.(\ref{Feq}) around a regular singular point $z= b_{i'}$ is a unit matrix, if and only if the monodromy matrix to Eq.(\ref{eq:Hg}) around a regular singular point $x= \pm \delta _{i'}$ is a unit matrix.
It is obvious that, if the monodromy matrix to Eq.(\ref{Feq}) around a regular singular point $z= b_{i'}$ is a unit matrix, then we have $r_{i'} \in \Zint _{\neq 0}$. In this paper we assume that $r_{i'} \in \Zint _{>0}$ for all $i'$.

\section{Integral representation and the Hermite-Krichever Ansatz} \label{sec:HK}

We introduce doubly-periodic functions to obtain an integral expression of solutions to Eq.(\ref{eq:H}) for the case $l_i \in \Zint _{\geq 0}$ $(i=0,1,2,3)$, $r_{i'} \in \Zint _{> 0}$ $(i'=1,\dots ,M)$ and the regular singular points $z=b_{i'}$ $(i'=1,\dots ,M)$ of Eq.(\ref{Feq}) are apparent.

\begin{prop} \label{prop:prod}
If $l_i \in \Zint _{\geq 0}$ $(i=0,1,2,3)$, $r_{i'} \in \Zint _{> 0}$ $(i'=1,\dots ,M)$ and regular singular points $z=b_{i'}$ $(i'=1,\dots ,M)$ of Eq.(\ref{Feq}) are apparent, then the equation
\begin{align}
& \left\{ \frac{d^3}{dx^3}-4\left( v(x) -E\right)\frac{d}{dx} -2\frac{dv(x)}{dx} \right\} \Xi (x)=0,
\label{prodDE}
\end{align}
has an even nonzero doubly-periodic solution that has the expansion
\begin{equation}
\Xi (x)=c_0+\sum_{i=0}^3 \sum_{j=0}^{l_i-1} b^{(i)}_j \wp (x+\omega_i)^{l_i-j} + \sum _{i'=1}^M  \sum_{j=0}^{r_{i'}-1} \frac{d^{(i')}_j}{(\wp (x)-\wp (\delta _{i'}))^{r_{i'} -j}}.
\label{Fx}
\end{equation}
\end{prop}
\begin{proof}
First, we show a lemma that is related to the monodromy of solutions to Eq. (\ref{eq:Hg}).
\begin{lemma} \label{prop:locmonod}
If $l_0, l_1, l_2, l_3  \in \Zint _{\geq 0}$, then the monodromy matrix of Eq.(\ref{eq:Hg}) around a point $x=n_1 \omega_1 + n_3 \omega_3$ $(n_1, n_3 \in \Zint)$ is a unit matrix.
\end{lemma} 
\begin{proof}
Due to periodicity, it is sufficient to consider the case $x=\omega _i$ $(i=0,1,2,3)$. We first deal with the case $i=1,2,3$.
The exponents at the singular point $x=\omega_i$ ($i=1,2,3$) are $-l_i$ and $l_{i}+1$. Because Eq.(\ref{eq:Hg}) is invariant under the transformation $x -\omega _i \rightarrow -(x-\omega _i) $ and the gap of the exponents at $x=\omega_i$ (i.e. $l_i+1-(-l_i)$) is odd, there exist solutions in the form $f_{i,1}(x)=(x-\omega _i)^{-l_i}(1+\sum_{j=1}^{\infty}a_j (x-\omega _i)^{2j})$ and $f_{i,2}(x)=(x-\omega _i)^{l_i+1}(1+\sum_{j=1}^{\infty}a'_j (x-\omega _i)^{2j})$.
Since the functions $f_{i,1}(x)$ and $f_{i,2}(x)$ form a basis for solutions to Eq.(\ref{eq:Hg}) and they are holomorphic around the point $x=\omega _i$, the monodromy matrix around $x=\omega _i$ is a unit matrix.
For the case $i=0$, the exponents at $x=0$ are $-l_0 - \sum_{i'=1}^M r_{i'}$ and $l_{0}+1 - \sum_{i'=1}^M r_{i'}$, and similarly it is shown that the monodromy matrix around the point $x=0$ is a unit matrix.
Hence we obtain the lemma.
\end{proof}
We continue the proof of Proposition \ref{prop:prod}.
Let $M_i$ $(i=1,3)$ be the transformations
obtained by the analytic continuation $x \rightarrow x+2\omega _i$.
It follows from double-periodicity of Eq.(\ref{eq:Hg}) that, if $f_g(x)$ is a solutions to Eq.(\ref{eq:Hg}), then $M_i f_g(x)$ $(i=1,3)$ is also a solution to Eq.(\ref{eq:Hg}).
From the assumption that regular singular points $z=b_{i'}$ are apparent for all ${i'}$, the monodromy matrix to Eq.(\ref{eq:Hg}) around a regular singular point $x= \pm \delta _{i'}$ is a unit matrix for all $i'$.
By combining with Lemma \ref{prop:locmonod}, it follows that all local monodromy matrices around any singular points are units. Hence the transformations $M_i$ do not depend on the choice of paths. From the fact that the fundamental group of the torus is commutative, we have $M_1 M_3=M_3 M_1$.
Recall that the operators $M_i$ act on the space of solutions to Eq.(\ref{eq:Hg}) for each $E$, which is two dimensional.
By the commutativity $M_1 M_3=M_3 M_1$, there exists a joint eigenvector $\tilde{\Lambda }_g (x)$ for the operators $M_1$ and $M_3$. 
It follows from Proposition \ref{prop:locmonod} and the apparency of singular points that the function $\tilde{\Lambda } _g (x)$ is single-valued and satisfies equations $(H_g-E-C_g) \tilde{\Lambda } _g(x)=0$, $M_1\tilde{\Lambda } _g(x)=\tilde{m}_1\tilde{\Lambda } _g(x)$ and $M_3\tilde{\Lambda } _g(x)=\tilde{m}_3\tilde{\Lambda } _g(x)$ for some $\tilde{m}_1,\tilde{m}_3 \in \Cplx \setminus \{0\}$.
By changing parity $x \leftrightarrow -x$, it follows immediately that $(H_g-E-C_g) \tilde{\Lambda } _g(-x)=0$, $M_1\tilde{\Lambda } _g(-x)=\tilde{m}_1^{-1}\tilde{\Lambda } _g(-x)$ and $M_3\tilde{\Lambda } _g(-x)=\tilde{m}_3^{-1}\tilde{\Lambda } _g(-x)$. Then the function $\tilde{\Lambda } _g(x)\tilde{\Lambda } _g(-x)$ is single-valued, even and doubly-periodic.
We set $\tilde{\Lambda } (x)= \tilde{\Lambda } _g(x) /\Psi _g (x) $. Then $\tilde{\Lambda } (x)$ and $\tilde{\Lambda } (-x)$ are solutions to Eq.(\ref{eq:H}).

Now consider the function $\Xi (x)=\tilde{\Lambda } _g(x)\tilde{\Lambda } _g(-x) /\Psi _g (x)^2$. Since the function $\Psi _g (x)^2$ is single-valued, even and doubly-periodic, the function $\Xi (x)$ is single-valued, even (i.e. $\Xi (x)=\Xi (-x)$), doubly-periodic (i.e. $\Xi (x+2\omega _1)=\Xi (x+2\omega _3)=\Xi (x)$), and satisfies the equation
\begin{align}
& \left\{ \frac{d^3}{dx^3}-4\left( v(x) -E\right)\frac{d}{dx} -2\frac{dv(x)}{dx} \right\} \Xi (x)=0
\nonumber
\end{align}
that the products of any pair of solutions to Eq.(\ref{eq:H}) satisfy.

Since the function $\Xi (x)$ is an even doubly-periodic function that satisfies the differential equation (\ref{prodDE}) and the exponents of Eq.(\ref{prodDE}) at $x=\omega _i$ $(i=0,\dots ,3)$ (resp. $x=\pm \delta _{i'}$ $(i'=1,\dots ,M)$) are $-2l_i,1,2l_i+2$ (resp. $-r _{i'} ,1,r _{i'} +2 $), it is written as a rational function of variable $\wp (x)$, and it admits the expansion as Eq.(\ref{Fx}) by considering exponents.
\end{proof}

The function $\Xi (x)$ is calculated by substituting Eq.(\ref{Fx}) into the differential equation (\ref{prodDE}) and solving simultaneous equations for the coefficients. We introduce an integral formula for a solution to the differential equation Eq.(\ref{eq:H}) in use of the function $\Xi (x)$.
Set
\begin{align}
& Q= \Xi (x)^2\left( E- v(x)\right) +\frac{1}{2}\Xi (x)\frac{d^2\Xi (x)}{dx^2}-\frac{1}{4}\left(\frac{d\Xi (x)}{dx} \right)^2. \label{const}
\end{align}
It follows from Eq.(\ref{prodDE}) that
\begin{align*}
\frac{dQ}{dx} = \frac{1}{2} \Xi (x) \left( 4\frac{d\Xi (x)}{dx} (E-v(x) )- 2 \Xi (x) \frac{dv(x)}{dx} + \frac{d^3\Xi (x)}{dx^3} \right) =0 .
\end{align*}
Hence the value $Q$ is independent of $x$.
\begin{prop} \label{prop:Linteg}
Let $\Xi (x)$ be the doubly-periodic function defined in Proposition \ref{prop:prod} and $Q$ be the value defined in Eq.(\ref{const}).
Then the function 
\begin{equation}
\Lambda ( x)=\sqrt{\Xi (x)}\exp \int \frac{ \sqrt{-Q}dx}{\Xi (x)},
\label{integ1}
\end{equation}
is a solution to the differential equation (\ref{eq:H}), and the function
\begin{equation}
\Lambda _g ( x)=\Psi _g (x) \sqrt{\Xi (x)}\exp \int \frac{ \sqrt{-Q}dx}{\Xi (x)},
\label{integ1g}
\end{equation}
is a solution to the differential equation (\ref{eq:Hg}).
\end{prop}
\begin{proof}
From Eqs.(\ref{integ1}, \ref{const}) we have
\begin{align}
& \frac{\Lambda ' (x)}{\Lambda (x)} =\frac{1}{2} \frac{\Xi ' (x)}{\Xi (x)} + \frac{\sqrt{-Q}}{\Xi (x)} , \label{eq:llp} \\
& \frac{\Lambda '' (x)}{\Lambda (x)} = \frac{1}{2} \frac{\Xi '' (x)}{\Xi (x)} - \frac{1}{4} \left( \frac{\Xi ' (x)}{\Xi (x)} \right)^2 - \frac{Q}{\Xi (x)^2} = v(x) -E.
\end{align}
Hence we have $-\frac{d^2}{dx^2} \Lambda (x) +v(x)\Lambda (x) =E\Lambda (x)$.
It follows from the equivalence of Eq.(\ref{eq:H}) and Eq.(\ref{eq:Hg}) that the function $\Lambda _g ( x)$ is a solution to Eq.(\ref{eq:Hg}). 
\end{proof}

\begin{prop} \label{prop:indep}
If $Q\neq 0$, then the functions $\Lambda (x) $ and $\Lambda (-x) $ are linearly independent and any solution to Eq.(\ref{eq:H}) is written as a linear combination of $\Lambda (x) $ and $\Lambda (-x) $.
\end{prop}
\begin{proof}
It follows from Eq.(\ref{eq:llp}) and the evenness of the function $\Xi (x)$ that 
\begin{equation}
\frac{\frac{d}{dx} \Lambda (-x)}{\Lambda (-x) }  = \frac{1}{2} \frac{\Xi ' (x)}{\Xi (x)} - \frac{\sqrt{-Q}}{\Xi (x)}.
\end{equation}
Hence we have
\begin{equation}
\Lambda (-x) \frac{d}{dx} \Lambda (x) -\Lambda (x) \frac{d}{dx} \Lambda (-x) = \Lambda (x) \Lambda (-x) \frac{2\sqrt{-Q}}{\Xi (x)}.
\label{eq:L-L}
\end{equation}
If $\Lambda (x) $ and $\Lambda (-x) $ are linearly dependent, then the l.h.s. of Eq.(\ref{eq:L-L}) must be zero; however, this is impossible because $Q \neq 0$. Hence the functions $\Lambda (x)$ and $\Lambda (-x)$ are linearly independent.
It follows from the invariance of Eq.(\ref{eq:H}) with respect to the transformation $x \leftrightarrow -x$ that $\Lambda (-x) $ is also a solution to Eq.(\ref{eq:H}).

Since solutions to Eq.(\ref{eq:H}) form a two-dimensional vector space and the functions $\Lambda (x)$ and $\Lambda (-x)$ are linearly independent, the functions $\Lambda (x)$ and $\Lambda (-x)$ form a basis of the space of solutions to Eq.(\ref{eq:H}), and any solution to Eq.(\ref{eq:H}) is written as a linear combination of $\Lambda (x) $ and $\Lambda (-x) $.
\end{proof}

It follows from Proposition \ref{prop:indep} that, if $Q\neq 0$, then the functions $\Lambda _g(x) $ and $\Lambda _g(-x) $ are linearly independent, and any solution to Eq.(\ref{eq:Hg}) is written as a linear combination of $\Lambda _g(x) $ and $\Lambda _g(-x) $.

From the formulae (\ref{integ1}, \ref{integ1g}) and the doubly-periodicity of the functions $\Xi (x)$ and $\Psi_g (x)^2$, we have 
\begin{align}
& \Lambda (x+2\omega _j)=\pm \Lambda (x) \exp \int _{0+ \varepsilon }^{2\omega _j+\varepsilon }\frac{\sqrt{-Q}dx}{\Xi (x)}, \quad (j=1,3), \label{eqn:Lam01} \\
& \Lambda _g (x+2\omega _j)=\pm \Lambda _g (x) \exp \int _{0+ \varepsilon }^{2\omega _j+\varepsilon }\frac{\sqrt{-Q}dx}{\Xi (x)}, \quad (j=1,3), \label{eqn:Lam01g}
\end{align}
with $\varepsilon $ a constant determined so as to avoid passing through the poles while integrating. The sign $\pm$ is determined by the analytic continuation of the function $\sqrt{\Xi (x)}$, and the integrations in Eqs.(\ref{eqn:Lam01}, \ref{eqn:Lam01g}) may depend on the choice of the path.
The function $\Lambda (x)$ may have branching points, althought the function $\Lambda _g (x)$ does not have branching points and is meromorphic on the complex plane, because $\Lambda _g (x)$ is a solution to Eq.(\ref{eq:Hg}) and any singularity of Eq.(\ref{eq:Hg}) is apparent.
It follows from Eq.(\ref{eqn:Lam01g}) that there exists $m_1, m_3 \in \Cplx $ such that 
\begin{equation}
\Lambda _g (x+2\omega _j) =\exp (\pi \sqrt{-1} m_j ) \Lambda _g (x), \quad (j=1,3).
\label{eqn:Lam010}
\end{equation}

We now show that a solution to Eq.(\ref{eq:H}) can be expressed in the form of the Hermite-Krichever Ansatz.
We set
\begin{equation}
\Phi _i(x,\alpha )= \frac{\sigma (x+\omega _i -\alpha ) }{ \sigma (x+\omega _i )} \exp (\zeta( \alpha )x), \quad \quad (i=0,1,2,3),
\label{Phii}
\end{equation}
where $\sigma (x)$ (resp. $\zeta (x)$) is the Weierstrass sigma (resp. zeta) function.
Then we have 
\begin{equation}
\left( \frac{d}{dx} \right) ^{k} \Phi _i(x+2\omega _{j} , \alpha ) = \exp (-2\eta _{j} \alpha +2\omega _{j} \zeta (\alpha )) \left( \frac{d}{dx} \right) ^{k} \Phi _i(x, \alpha )
\label{ddxPhiperiod}
\end{equation}
for $i=0,1,2,3$, $j=1,2,3$ and $k \in \Zint _{\geq 0}$, where $\eta _j =\zeta (\omega _j)$ $(j=1,2,3)$.
\begin{thm} \label{thm:alpha}
Set $\tilde{l} _0 = l_0 +\sum _{i'=1}^M r_{i'}$ and $\tilde{l}_i =l_i$ $(i=1,2,3)$.
The function $\Lambda _g (x)$ in Eq.(\ref{integ1g}) is expressed as
\begin{align}
& \Lambda _g (x) = \exp \left( \kappa x \right) \left( \sum _{i=0}^3 \sum_{j=0}^{\tilde{l}_i-1} \tilde{b} ^{(i)}_j \left( \frac{d}{dx} \right) ^{j} \Phi _i(x, \alpha ) \right)
\label{Lalpha}
\end{align}
for some values $\alpha $, $\kappa$ and $\tilde{b} ^{(i)}_j$ $(i=0,\dots ,3, \: j= 0,\dots ,\tilde{l}_i-1)$, or
\begin{align}
& \Lambda _g (x) = \exp \left( \bar{\kappa } x \right) \left( \bar{c} +\sum _{i=0}^3 \sum_{j=0}^{\tilde{l}_i-2} \bar{b} ^{(i)}_j \left( \frac{d}{dx} \right) ^{j} \wp (x+\omega _i) +\sum_{i=1}^3 \bar{c}_i \frac{\wp '(x)}{\wp (x)-e_i} \right)
\label{Lalpha0}
\end{align}
for some values $\bar{\kappa }$, $\bar{c}$, $\bar{c}_i$ $(i=1,2,3)$ and $\bar{b} ^{(i)}_j$ $(i=0,\dots ,3, \: j= 0,\dots ,\tilde{l}_i-2)$.

If the function $\Lambda _g (x)$ is expressed as Eq.(\ref{Lalpha}), then
\begin{align}
& \Lambda _g (x+2\omega _j) = \exp (-2\eta _j \alpha +2\omega _j \zeta (\alpha ) +2 \kappa \omega _j ) \Lambda _g (x) , \quad  (j=1,3), \label{ellint} 
\end{align}
else
\begin{align} 
& \Lambda _g (x+2\omega _j) = \exp (2 \bar{\kappa } \omega _j ) \Lambda _g (x) , \quad  (j=1,3). \label{ellint0} 
\end{align}		
\end{thm}

\begin{proof}
Set
\begin{align}
& \alpha  = -m_1 \omega _3 +m_3 \omega _1 \label{al},
\end{align}
where $m_1$ and $m_3 $ are determined in Eq.(\ref{eqn:Lam010}).

If $\alpha \not \equiv 0$ $($mod $2\omega_1 \Zint \oplus 2\omega_3 \Zint)$, then we set
\begin{align}
 & \kappa = \zeta (m_1 \omega _3 -m_3 \omega _1 ) -m_1 \eta _3 +m_3 \eta _1 . \label{kapp}
\end{align}
It follows from Legendre's relation $\eta _1 \omega _3 - \eta_3 \omega _1 =\pi \sqrt{-1} /2$ that 
\begin{align}
& \exp (\kappa  (x+2\omega _{j}) ) \left( \frac{d}{dx} \right) ^{k} \Phi _i(x+2\omega _{j} , \alpha ) \label{eq:periodj'} \\
& = \exp (-2\eta _{j} \alpha +2\omega _{j} (\zeta (\alpha ) + \kappa )) \exp (\kappa x )\left( \frac{d}{dx} \right) ^{k} \Phi _i(x, \alpha ) \nonumber \\
& = \exp (2m_1 (\eta _{j} \omega _3- \eta _3\omega _j) + 2m_3 (\eta _{1} \omega _j- \eta _j \omega _1)) \exp (\kappa x )\left( \frac{d}{dx} \right) ^{k} \Phi _i(x, \alpha ) \nonumber \\
& = \exp (\pi \sqrt{-1} m_{j}) \exp (\kappa x )\left( \frac{d}{dx} \right) ^{k} \Phi _i(x, \alpha ) \nonumber
\end{align}
for $i=0,1,2,3$, $j=1,3$ and $k \in \Zint _{\geq 0}$.
Hence the function $\Lambda _g (x)$ and the functions $\exp (\kappa x )\left( \frac{d}{dx} \right) ^{k} \Phi _i(x, \alpha )$ have the same periodicity with respect to periods $(2\omega_1, 2\omega _3)$.
Since the meromorphic function $\Lambda _g (x)$ satisfies Eq.(\ref{eq:Hg}), it is holomorphic except for $\Zint \omega _1 \oplus \Zint \omega _3 $ and has a pole of degree $\tilde{l} _i$ or zero of degree $\tilde{l} _i +1$ at $x= \omega _i$ $(i=0,1,2,3)$.
The function $\exp(\kappa x) \left( \frac{d}{dx} \right) ^{k} \Phi _i(x , \alpha ) $ has a pole of degree $k+1$ at $x=\omega _i$.
By subtracting the functions $\exp(\kappa x) \left( \frac{d}{dx} \right) ^{k} \Phi _i(x , \alpha ) $ from the function $\Lambda _g (x)$ to erase the poles, we obtain a holomorphic function that has the same periods as $\Phi _0(x , \alpha ) $, and must be zero.
Hence we obtain the expression (\ref{Lalpha}). 
The periodicity (see Eq.(\ref{ellint})) follows from Eq.(\ref{eq:periodj'}).

If $\alpha \equiv 0$ $($mod $2\omega_1 \Zint \oplus 2\omega_3 \Zint)$ (i.e. $m_1 \omega _3 \equiv m_3 \omega _1$ $($mod $2\omega_1 \Zint \oplus 2\omega_3 \Zint)$), then we set
\begin{align}
 & \bar{\kappa }= -m_1 \eta _3 +m_3 \eta _1 .
\end{align}
The function $\Lambda _g (x)$ and the function $\exp (\bar{\kappa }x )$ have the same periodicity with respect to periods $(2\omega_1, 2\omega _3)$.
Hence the function $\Lambda _g (x) \exp (-\bar{\kappa }x )$ is doubly periodic, and we obtain the expression (\ref{Lalpha0}) by considering the poles. Periodicity (see Eq.(\ref{ellint0})) follows immediately.
\end{proof}

We investigate the situation that Eq.(\ref{eq:Hg}) has a non-zero solution of an elliptic function.
Let ${\mathcal F}_{\epsilon _1 , \epsilon _3 }$ and $\tilde{\mathcal F} _{\epsilon _1 , \epsilon _3 }$ $(\epsilon _1 , \epsilon _3 \in \{ \pm 1 \})$ be the spaces defined by
\begin{align}
& {\mathcal F} _{\epsilon _1 , \epsilon _3 }=\{ f(x) \mbox{: meromorphic }| f(x+2\omega_1)= \epsilon _1 f(x), \; f(x+2\omega_3)= \epsilon _3 f(x) \} , \\
& \tilde{\mathcal F} _{\epsilon _1 , \epsilon _3 }=\left\{ f(x) \: \left| 
\begin{array}{l}
f(x) \Psi _g (x) \in {\mathcal F}_{\epsilon _1 , \epsilon _3 } \mbox{ and holomorphic}\\
\mbox{except for } \Zint \omega _1 \oplus \Zint \omega _3 \mbox{, and the degree of the pole at }\\
 x=\omega _i \mbox{ is no more than } \left\{ 
\begin{array}{ll}
l_i, &   i=1,2,3, \\
l_0+\sum _{i'=1}^M r_{i'}, & i=0.
\end{array}
\right.
\end{array}
\right. \right\} ,
\end{align}
where $(2\omega_1, 2\omega_3)$ are basic periods of elliptic functions. Then $\tilde{\mathcal F} _{\epsilon _1 , \epsilon _3 }$ is a finite-dimensional vector space. 
Note that, if a solution $f(x)$ to Eq.(\ref{eq:H}) satisfies the condition $f(x+2\omega _1) \Psi _g (x +2\omega _1) =\epsilon _1 f(x) \Psi _g (x) $ and $f(x+2\omega _3) \Psi _g (x +2\omega _3) =\epsilon _3 f(x) \Psi _g (x) $ for some $\epsilon _1 , \epsilon _3 \in \{ \pm 1 \}$, then we have $f(x) \in {\mathcal F} _{\epsilon _1 , \epsilon _3 }$, because the position of the poles and their degree are restricted by the differential equation.

\begin{prop} \label{prop:disttwoch} 
Assume that Eq.(\ref{eq:H}) has a non-zero solution in the space $\tilde{\mathcal F} _{\epsilon _1 , \epsilon _3 }$ for some $\epsilon _1 , \epsilon _3 \in \{ \pm 1 \}$. Then the signs $(\epsilon _1 , \epsilon _3 )$ are determined uniquely for each $E$, $\tilde{s}_{i'}$ $(i'=1, \dots ,M)$ etc.
\end{prop}
\begin{proof}
Assume that Eq.(\ref{eq:H}) has a non-zero solution in both the spaces $\tilde{\mathcal F} _{\epsilon _1 , \epsilon _3 }$ and $\tilde{\mathcal F} _{\epsilon '_1 , \epsilon '_3 }$. Let $f_1 (x)$ (resp. $f_2(x)$) be the solution to the differential equation (\ref{eq:H}) in the space $\tilde{\mathcal F} _{\epsilon _1 , \epsilon _3 }$ (resp. the space $\tilde{\mathcal F} _{\epsilon '_1 , \epsilon '_3 }$). 
Then periodicity of the function $f_1 (x)\Psi _g (x)$ and $f_2(x)\Psi _g (x)$ is different, more precisely there exists $j \in \{ 1,3\}$ such that 
\begin{equation}
\left\{ 
\begin{array}{ll}
f_1 (x +2\omega _j) \Psi _g (x +2\omega _j)= \pm f_1 (x) \Psi _g (x),\\
f_2 (x +2\omega _j) \Psi _g (x +2\omega _j)= \mp f_2 (x) \Psi _g (x).
\end{array}
\right.
\label{f1f2antp}
\end{equation}
Then the functions $f_1(x)$ and $f_2(x)$ are linearly independent.
Since the functions $f_1(x)$ and $f_2(x)$ satisfy Eq.(\ref{eq:H}), we have $\frac{d}{dx} \left( f_2(x) f'_1(x) -f_1(x) f'_2(x) \right) = f_2(x) f''_1(x) -  f_1(x) f''_2(x)= 0$.
Therefore $f_2(x) f'_1 (x) - f_1(x) f'_2 (x) =C$ for constants $C$, and $C$ is non-zero, which follows from linear independence. By Eq.(\ref{f1f2antp}), the function $(f_2(x) f'_1 (x) - f_1(x) f'_2 (x) )\Psi _g (x)^2$ is anti-periodic with respect to the period $2\omega _j$, but it contradicts to $C\neq 0$.
Hence, we proved that Eq.(\ref{eq:H}) does not have a non-zero solution in both the spaces $\tilde{\mathcal F} _{\epsilon _1 , \epsilon _3 }$ and $\tilde{\mathcal F} _{\epsilon '_1 , \epsilon '_3 }$.
\end{proof}

\begin{prop} \label{prop:Q0F}
If $Q=0$, then we have $\Lambda (x) \in \tilde{\mathcal F} _{\epsilon _1 , \epsilon _3 }$ for some $\epsilon _1 , \epsilon _3 \in \{ \pm 1 \}$.
\end{prop}
\begin{proof}
It follows from Eq.(\ref{integ1}) and the double-periodicity of the function $\Xi (x)\Psi _g (x)^2$ that
\begin{equation}
(\Lambda (x+2\omega _j) \Psi _g (x+2\omega _j))^2 = (\Lambda (x)\Psi _g (x) )^2 = \Xi (x)\Psi _g (x)^2,
\end{equation}
for $j=1,3$. Hence $\Lambda (x+2\omega _j) \Psi _g (x+2\omega _j)= \pm \Lambda (x)\Psi _g (x) $ $(j=1,3)$
and we have $\Lambda (x) \in  \tilde{\mathcal F} _{\epsilon _1 , \epsilon _3 }$ for some $\epsilon _1 , \epsilon _3 \in \{ \pm 1 \}$.
\end{proof}

It follows from Proposition \ref{prop:prod} that the dimension of the space of solutions to Eq.(\ref{prodDE}), which are even doubly-periodic, is no less than one.
Since the exponents of Eq.(\ref{prodDE}) at $x=0$ are $-2l_0$, $1$ and $2l_0 +2$, the dimension of the space of even solutions to Eq.(\ref{prodDE}) is at most two.
Hence, the dimension of the space of solutions to Eq.(\ref{prodDE}), which are even doubly-periodic, is one or two.

\begin{prop} \label{prop:twoF}
Assume that the dimension of the space of solutions to Eq.(\ref{prodDE}), which are even doubly-periodic, is two. The all solutions to Eq.(\ref{eq:H}) are contained in the space $\tilde{\mathcal F} _{\epsilon _1 , \epsilon _3 }$ for some $\epsilon _1 , \epsilon _3 \in \{ \pm 1 \}$.
\end{prop}
\begin{proof}
Since the differential equation (\ref{eq:H}) is invariant under the change of parity $x \leftrightarrow -x$, a basis of the solutions to Eq.(\ref{eq:H}) is taken as $f_e (x)$ and $f_o (x)$ such that $f_ e(x)$ (resp. $f_o (x)$) satisfies $f_ e(-x)=f_e (x)$ (resp. $f_ o(-x)=-f_o (x)$).
Then the functions $f_ e(x)^2$ and $f_ o(x)^2$ are even and they are solutions to Eq.(\ref{prodDE}). Since the dimension of the space of even solutions to Eq.(\ref{prodDE}) is at most two, and the dimension of the space of solutions to Eq.(\ref{prodDE}), which are even doubly-periodic, is two, the even functions $f_e(x)^2$ and $f_o(x)^2$ must be doubly-periodic. Hence $(f_e(x+2\omega _j) \Psi _g (x+2\omega _j))^2 =(f_e(x) \Psi _g (x))^2$ $(j=1,3)$ and it follows that $f_e(x+2\omega _j) \Psi _g (x+2\omega _j)= \pm f_e(x)\Psi _g (x)$ $(j=1,3)$. Therefore we have $f_e (x) \in \tilde{\mathcal F} _{\epsilon _1 , \epsilon _3 }$ for some $\epsilon _1 , \epsilon _3 \in \{ \pm 1 \}$. Similarly we have $f_ o(x) \in \tilde{\mathcal F} _{\epsilon '_1 , \epsilon '_3 }$ for some $\epsilon '_1 , \epsilon '_3 \in \{ \pm 1 \}$, and it follows from Proposition \ref{prop:disttwoch} that $\epsilon '_j =\epsilon _j$ $(j=1,3)$. Since $f_e (x)$ and $f_o (x)$ are a basis of solutions to Eq.(\ref{eq:H}), all solutions to Eq.(\ref{eq:H}) are contained in the space $\tilde{\mathcal F} _{\epsilon _1 , \epsilon _3 }$.
\end{proof}

\begin{prop} \label{prop:Xionedim}
If $M=0$ or ($M=1$ and $r_1 =1$), then the dimension of the space of solutions to Eq.(\ref{prodDE}), which are even doubly-periodic, is one.
\end{prop}
\begin{proof}
Assume that the dimension of the space of solutions to Eq.(\ref{prodDE}), which are even doubly-periodic, is two.
From Proposition \ref{prop:twoF}, all solutions to Eq.(\ref{eq:H}) are contained in the space $\tilde{\mathcal F} _{\epsilon _1 , \epsilon _3 }$ for some $\epsilon _1 , \epsilon _3 \in \{ \pm 1 \}$. Since the differential equation (\ref{eq:Hg}) is invariant under the change of parity $x \leftrightarrow -x$, a basis of the solutions to Eq.(\ref{eq:Hg}) is taken as $f_1(x)$ and $f_2 (x)$ such that $f_1 (x)$ (resp. $f_2 (x)$) is even (resp. odd) function.
From the assumption that $l_i \in \Zint $ $(i=0,1,2,3)$ and that regular singular points $b_{i'}$ are apparent ($i'=1,\dots ,M$), the functions $f_1 (x)$ and $f_2 (x)$ are meromorphic.
Since the function $f_1 (x)$ (resp. $f_2 (x)$) satisfies Eq.(\ref{eq:Hg}), it does not have poles except for $\Zint \omega _1 \oplus \Zint \omega _3$.
Hence the function $f_1 (x)$ admits the expression $f_1 (x)= \wp _1 (x) ^{\tilde{\beta }_1}\wp _2 (x) ^{\tilde{\beta }_2}  \wp _3 (x) ^{\tilde{\beta }_3} (P^{(1)}(\wp (x))+\wp '(x) P^{(2)}(\wp (x)))$, where $\wp _i (x)$ $(i=1,2,3)$ are co-$\wp $ functions and $P^{(1)}(z)$, $P^{(2)} (z)$ are polynomials in $z$. Since the function $f_1 (x) $ is even, we have $P^{(1)} (z)=0$ or $P^{(2)} (z)=0$. By combining with the relation $\wp '(z) = -2\wp _1 (z) \wp _2 (z) \wp _3 (z) $, the function $f_1 (x)$ is expressed as 
\begin{equation}
f_1 (x)= \wp _1 (x) ^{\beta _1}\wp _2 (x) ^{\beta _2}  \wp _3 (x) ^{\beta _3} P_1 (\wp (x)),
\end{equation}
where $P_1(z)$ is a polynomial in $z$.
Because the exponents of Eq.(\ref{eq:Hg}) at $x=\omega _i$ $(i=1,2,3)$ are $-l_i$ and $l_{i} +1$, we have $\beta _i \in \{ -l_i ,l _{i} +1 \}$ $(i=1,2,3)$.
Similarly the function $f_2 (x)$ is expressed as 
\begin{equation}
f_2(x)= \wp _1 (x) ^{\beta '_1}\wp _2 (x) ^{\beta '_2}  \wp _3 (x) ^{\beta '_3} P_2 (\wp (x)),
\end{equation}
where $P_2(z)$ is a polynomial in $z$ and $\beta ' _i \in \{ -l_i ,l _{i} +1 \}$.

Since the functions $ \wp _i (x) $ $(i=1,2,3)$ are odd and the parity of functions $f_1 (x)$ and $f_2 (x)$ is different, we have $\beta _1 +\beta _2 +\beta _3 \not \equiv \beta '_1 +\beta '_2 +\beta '_3 $ (mod $2$). Since $f_j ( x +2\omega _1) = (-1) ^{\beta _2 +\beta _3 } f_j (x)$, $ f_j ( x +2\omega _3) = (-1) ^{\beta _1 +\beta _2 } f_j (x)$ $(j=1,2)$, we have $\beta _2 +\beta _3 \equiv \beta '_2 +\beta '_3 $ (mod $2$) and $\beta _1 +\beta _2 \equiv \beta '_1 +\beta '_2 $ (mod $2$). Hence we have $\beta _i \not \equiv \beta '_i $ (mod $2$) for $i=1,2,3$. Therefore $(\beta _i , \beta '_i ) = (-l_i , l _{i} +1)$ or $(\beta _i , \beta '_i ) = (l _{i} +1, -l_i)$ for each $i \in \{ 1,2,3 \}$.
Let $\beta _0$ (resp. $\beta '_0$) be the exponent of the function $f_1 (x)$ (resp. $f_2 (x)$) at $x=0$. Since the parity of functions $f_1 (x)$ and $f_2 (x)$ is different and the exponents of Eq.(\ref{eq:Hg}) at $x=0$ are $-l_0- \sum_{i'=1}^M r_{i'}$ and $l _{0} +1- \sum_{i'=1}^M r_{i'}$, we have $(\beta _0 , \beta '_0 ) = (-l_0 - \sum_{i'=1}^M r_{i'}, l _0 +1- \sum_{i'=1}^M r_{i'})$ or $(\beta _0, \beta '_0 ) = ( l _0 +1- \sum_{i'=1}^M r_{i'} ,-l_0 - \sum_{i'=1}^M r_{i'})$.

Since the function $f_1 (x)$ is doubly-periodic with periods $(4\omega _1 , 4\omega _3)$, the sum of degrees of zeros of $f_1 (x)$ on the basic domain is equal to the sum of degrees of poles of $f_1 (x)$. Since the function $f_1 (x)$ does not have poles except for $\Zint \omega _1 \oplus \Zint \omega _3$, we have $\sum _{i=0}^3 \beta _i \leq 0$. Similarly we have $\sum _{i=0}^3 \beta '_i \leq 0$.
Hence $0 \geq \sum _{i=0}^3 (\beta _i +\beta '_i) = 4-2 \sum_{i'=1}^M r_{i'}$.
Therefore we have $\sum_{i'=1}^M r_{i'} \geq 2$. 

Thus we obtain that, if $M=0$ or ($M=1$ and $r_1 =1$), then the dimension of the space of solutions to Eq.(\ref{prodDE}), which are even doubly-periodic, is one.
\end{proof}

Note that the case $M=0$ corresponds to Heun's equation, and the case $M=1$ and $r_1 =1$ is related with the sixth Painlev\'e equation.

\begin{exa}
Let us consider the following differential equation:
\begin{equation}
\left\{ -\left( \frac{d}{dx} \right) ^2 +\left( \frac{\wp ' (x)}{\wp (x) + \sqrt{\frac{g_2}{12}}}+\frac{\wp ' (x)}{\wp (x) - \sqrt{\frac{g_2}{12}}} \right) \frac{d}{dx} \right\} f(x)=0
\label{deq:twodim}
\end{equation}
This equation corresponds to the case $l_0=1$, $l_1=l_2=l_3=0$, $M=2$ and $r_1=r_2=1$, if $g_2 \neq 0$.
From the relation
\begin{equation}
 \frac{\wp ' (x)}{\wp (x) + \sqrt{\frac{g_2}{12}}}+\frac{\wp ' (x)}{\wp (x) - \sqrt{\frac{g_2}{12}}} = \frac{\wp ''' (x)}{\wp '' (x)},
\end{equation}
a basis of the solutions to Eq.(\ref{deq:twodim}) is $1$, $\wp '(x)$.
The dimension of the solutions to Eq.(\ref{prodDE}), which are even doubly-periodic, is two, and a basis of the solutions to Eq.(\ref{prodDE}) is written as $1/\wp ''(x)$, $\wp '(x)^2/\wp ''(x)$, $\wp' (x)/\wp ''(x)$.
\end{exa}

\begin{prop} \label{prop:char}
Assume that the dimension of the space of the solutions to Eq.(\ref{prodDE}), which are even doubly-periodic, is one. 
Let $c_0$, $b^{(i)}_j$ and $d^{(i')}_j$ be constants defined in Eq.(\ref{Fx}).\\
(i) If there exists a non-zero solution to Eq.(\ref{eq:H}) in the space $\tilde{\mathcal F} _{\epsilon _1 , \epsilon _3 }$ for some $\epsilon _1 , \epsilon _3 \in \{ \pm 1 \}$, then we have $Q=0$.\\
(ii) If $Q \neq 0$ and $l_i \neq 0$, then $b^{(i)}_0 \neq 0$.\\
(iii) If $Q \neq 0$ and $l_0 =0$, then $\Xi (0)\neq 0$. In particular, if $Q \neq 0$ and $l_0=l_1=l_2=l_3 =0$, then $c_0 \neq 0$.
\end{prop}
\begin{proof}
First we prove (i). Suppose that there exists a non-zero solution to Eq.(\ref{eq:H}) in the space $\tilde{\mathcal F} _{\epsilon _1 , \epsilon _3 }$ and $Q \neq 0$. From the condition $Q \neq 0$, the functions $\Lambda (x) $ and $\Lambda (-x) $ form the basis of the space of the solutions to the differential equation (\ref{eq:H}).
Since there is a non-zero solution to Eq.(\ref{eq:H}) in the space $\tilde{\mathcal F} _{\epsilon _1 , \epsilon _3 }$, there exist constants $(C_1, C_2) \neq (0,0)$ such that $C_1 \Lambda (x) +C_2 \Lambda (-x) \in \tilde{\mathcal F} _{\epsilon _1 , \epsilon _3 }$.
By shifting $x \rightarrow x +2\omega _i$ ($i=1,3$), it follows from Eq.(\ref{eqn:Lam010}) that
\begin{align}
& \quad (C_1 \Lambda (x +2\omega _i ) +C_2 \Lambda (-(x+2\omega _i) )) \Psi _g (x+2\omega _i) \\
& = C_1 \Lambda (x +2\omega _i ) \Psi _g (x+2\omega _i) \pm C_2 \Lambda (-x-2\omega _i) ) \Psi _g (-x-2\omega _i) \nonumber \\
& = C_1 \exp (\pi \sqrt{-1} m_i ) \Lambda (x) \Psi _g (x) \pm C_2 \exp (-\pi \sqrt{-1} m_i ) \Lambda (-x) \Psi _g (-x) \nonumber \\
& = (C_1 \exp (\pi \sqrt{-1} m_i ) \Lambda (x) +C_2 \exp (-\pi \sqrt{-1} m_i ) \Lambda (-x) ) \Psi _g (x) , \nonumber
\end{align}
where the sign $\pm $ is determined by the branching of the function $\Psi _g (x)$, and the function $C_1 \exp (\pi \sqrt{-1} m_i ) \Lambda (x) +C_2 \exp (-\pi \sqrt{-1} m_i ) \Lambda (-x)$ also satisfies Eq.(\ref{eq:H}).
On the other hand, it follows from the definition of the space $\tilde{\mathcal F} _{\epsilon _1 , \epsilon _3 }$ that $(C_1 \Lambda (x +2\omega _i ) +C_2 \Lambda (-(x+2\omega _i) )) \Psi _g (x+2\omega _i)= \pm (C_1 \Lambda (x ) +C_2 \Lambda (-x ))\Psi _g (x)$ for signs $\pm $. By comparing two expressions, we have $\exp (\pi \sqrt{-1} m_i ) \in \{ \pm 1 \}$ ($i=1,3$) and the periodicities of the functions $\Lambda (x ) \Psi _g (x)$  and $(C_1 \Lambda (x ) +C_2 \Lambda (-x ))\Psi _g (x)$ coincide. Thus $ \Lambda (x ), \: \Lambda (-x ) \in \tilde{\mathcal F} _{\epsilon _1 , \epsilon _3 }$. The functions $\Lambda (x )^2$ and $\Lambda (-x )^2$ are even doubly-periodic function and satisfy Eq.(\ref{prodDE}), because they are the products of a pair of solutions to Eq.(\ref{eq:H}). Hence the dimension of the space of solutions to Eq.(\ref{eq:H}), which are even doubly-periodic, is no less than two, and contradict the assumption of the proposition.
Therefore the supposition $Q\neq 0$ is false, and we obtain (i).

Next we show (ii). Assume that $l_i \neq 0$. Since the exponents of Eq.(\ref{eq:H}) at $x=\omega_{i}$ are $-l_i$ or $l_i +1$, the function $\Lambda (x)$ has a pole of degree $l_i$ or a zero of degree $l_i +1$ at $x=\omega_{i}$. It follows from the periodicity (see Eq.(\ref{eqn:Lam010})) that, if the function $\Lambda (x)$ has a zero at $x=\omega_{i}$, then $\Lambda (x)$ has also a zero at $x=-\omega_{i}$.
Hence the function $\Lambda (-x)$ has a zero at $x=\omega_{i}$. From the assumption $Q \neq 0$, any solution to Eq.(\ref{eq:H}) is written as a linear combination of functions $\Lambda (x)$ and $\Lambda (-x)$. But it contradicts that one of the exponents at $x=\omega_{i}$ is $-l_i$. Hence the function $\Lambda (x)$ has a pole of degree $l_i$ and $b^{(i)}_0 \neq 0$.

(iii) is proved similarly by showing that the function $\Lambda (x)$ does not have zero at $x=0$.
\end{proof}

By combining Propositions \ref{prop:Q0F} and \ref{prop:char} (i) we obtain the following proposition:
\begin{prop} \label{prop:zeros}
Assume that the dimension of the space of solutions to Eq.(\ref{prodDE}), which are even doubly-periodic, is one.
Then the condition $Q=0$ is equivalent to that there exists a non-zero solution to Eq.(\ref{eq:H}) in the space $\tilde{\mathcal F} _{\epsilon _1 , \epsilon _3 }$ for some $\epsilon _1 , \epsilon _3 \in \{ \pm 1 \}$.
\end{prop}

We show that the function $\Lambda (x) $ admits an expression of the Bethe Ansatz type.
\begin{prop} \label{prop:BA}
Set $l=\sum _{i=0}^3 l_i +\sum _{i'=1}^M r_{i'} $, $\tilde{l} _0 = l_0 +\sum _{i'=1}^M r_{i'}$ and $\tilde{l}_i =l_i$ $(i=1,2,3)$. Assume that  $Q \neq 0$ and the dimension of the space of the solutions to Eq.(\ref{prodDE}), which are even doubly-periodic, is one.\\
(i) The function $\Lambda (x) $ in Eq.(\ref{integ1}) is expressed as 
\begin{align}
& \Lambda (x) = \frac{C_0 \prod_{j=1}^l \sigma(x-t_j)}{\Psi _g (x) \sigma(x)^{\tilde{l}_0}\sigma_1(x)^{\tilde{l}_1}\sigma_2(x)^{\tilde{l}_2}\sigma_3(x)^{\tilde{l}_3}}\exp \left(cx \right), 
\label{eq:tilL}
\end{align}
for some $t_1, \dots , t_l$, $c$ and $C_0 (\neq 0)$, where $\sigma _i (x)$ $(i=1,2,3)$ are co-sigma functions.\\
(ii) $t_j + t_{j'} \not \equiv 0$ $($mod $2\omega_1 \Zint \oplus 2\omega_3 \Zint)$ for all $j,j'$.\\
(iii) If $t_j \not \equiv \pm \delta _{i'}$ $($mod $2\omega_1 \Zint \oplus 2\omega_3 \Zint)$ for all $i' \in \{1,\dots ,M\}$, then we have $t_j \not \equiv t_{j'}$ $($mod $2\omega_1 \Zint \oplus 2\omega_3 \Zint)$ for all $j' (\neq j)$.\\
(iv) If $t_j \equiv \pm \delta _{i'}$ $($mod $2\omega_1 \Zint \oplus 2\omega_3 \Zint)$, then $\# \{ j' \: | \: t_j \equiv t_{j'}$ $($mod $2\omega_1 \Zint \oplus 2\omega_3 \Zint) \} = r_{i'}+1$.
(v) If $l_0 \neq 0$ (resp. $l_0=0$), then we have $c=\sum_{i=1}^l \zeta(t_j)$ (resp. $c=\sum_{i=1}^l \zeta(t_j)+ \frac{\sqrt{-Q}}{\Xi (0)}$). (Note that it follows from Proposition \ref{prop:char} that $\frac{\sqrt{-Q}}{\Xi (0)}$ is finite.)\\
(vi) Set $z= \wp (x)$ and $z_j=\wp (t_j)$. Then 
\begin{equation}
\left. \frac{d\Xi (x)}{dz}\right| _{z=z_j}=\frac{2\sqrt{-Q}}{\wp'(t_j)}.
\label{signpptj}
\end{equation}
\end{prop}
\begin{proof}
Let $\alpha $ be the value defined in Eq.(\ref{al}).
First, we consider the case $\alpha \not \equiv 0$ $($mod $2\omega_1 \Zint \oplus 2\omega_3 \Zint)$. Let $\kappa $ be the value defined in Eq.(\ref{kapp}).
Then the function $\Lambda _g (x) / \left( \exp (\kappa x ) \Phi _0(x, \alpha ) \right) $ is meromorphic and doubly-periodic. Hence there exists $a_1 ,\dots ,a_{l'}$, $b_1 ,\dots ,b_{l'}$ such that $a_1 + \dots + a_{l'} = b_1 +\dots + b_{l'}$ and 
$$
\Lambda _g (x) / \left( \exp (\kappa x ) \Phi _0(x, \alpha ) \right)  = \frac{\prod_{i=1}^{l'} \sigma (x- a_i )}{\prod_{i=1}^{l'} \sigma (x- b_i )}.
$$
For the case $\alpha \equiv 0$ $($mod $2\omega_1 \Zint \oplus 2\omega_3 \Zint)$ the function $\Lambda _g (x) / \exp (\bar{\kappa }x)$ is similarly expressed as
$$
\Lambda _g (x) / \exp (\bar{\kappa }x) = \frac{\prod_{i=1}^{l'} \sigma (x- a_i )}{\prod_{i=1}^{l'} \sigma (x- b_i )}.
$$

Since the function $\Lambda _g (x)$ satisfies Eq.(\ref{eq:Hg}), it does not have poles except for $\omega_1 \Zint \oplus \omega_3 \Zint$. From Proposition \ref{prop:char} (ii), it has poles at $x=\omega _i$ of degree $\tilde{l}_i$. Hence we have the expression
\begin{align}
& \Lambda _g (x) = \frac{C_0 \prod_{j=1}^l \sigma(x-t_j)}{\sigma(x)^{\tilde{l}_0}\sigma_1(x)^{\tilde{l}_1}\sigma_2(x)^{\tilde{l}_2}\sigma_3(x)^{\tilde{l}_3}}\exp \left(cx \right), 
\label{eq:tilL0}
\end{align}
for some $t_1, \dots , t_l$, $c$ and $C_0 (\neq 0)$ such that $t_j \not \equiv 0$ $($mod $\omega_1 \Zint \oplus \omega_3 \Zint)$.
Therefore we obtain (i) and that $2t_j \not \equiv 0$ $($mod $2\omega_1 \Zint \oplus 2\omega_3 \Zint)$.

Suppose that $t_j + t_{j'} \equiv 0$ $($mod $2\omega_1 \Zint \oplus 2\omega_3 \Zint)$ for some $j $ and $j'$,
From Eq.(\ref{eq:tilL}) and $-t_j \equiv t_{j'} $ $($mod $2\omega_1 \Zint \oplus 2\omega_3 \Zint)$, we have $\Lambda _g (t_j) = \Lambda  _g (-t_{j}) =0$. Since $Q\neq 0$, all solutions to Eq.(\ref{eq:Hg}) are written as linear combinations of $\Lambda  _g (x)$ and $\Lambda  _g (-x) $. Hence $t_j$ is a zero for all solutions to Eq.(\ref{eq:Hg}), but they contradict that one of the exponents at $x=t_j$ is zero. Therefore we obtain (ii).

If $t_j \not \equiv \pm \delta _{i'}, \omega _i$ $($mod $2\omega_1 \Zint \oplus 2\omega_3 \Zint)$ for all $i$ and $i'$, then the exponents of Eq.(\ref{eq:Hg}) at $x=t_j$ are $0$ and $1$. Hence $x=t_j$ is a zero of $\Lambda _g (x)$ of degree one.
Incidentally, the exponents of Eq.(\ref{eq:Hg}) at $x=\pm \delta _{i'}$ are $0$ and $r_{i'}+1$. Hence, if $t_j \equiv \pm \delta _{i'}$ $($mod $2\omega_1 \Zint \oplus 2\omega_3 \Zint)$, then $x=t_j$ is a zero of $\Lambda _g (x)$ of degree $r_{i'}+1$. Thus we obtain (iii) and (iv).

It follows from Eq.(\ref{eq:tilL}) that
\begin{equation}
\frac{\Lambda ' (x)}{\Lambda (x)} = c- \tilde{l}_0 \frac{\sigma '(x)}{\sigma (x)} -\sum _{i=1}^3 \tilde{l}_i \frac{\sigma ' _i (x)}{\sigma _i (x)}+\sum_{j=1}^l \frac{\sigma '(x-t_j)}{\sigma (x-t_j)} -\sum _{i' =1}^M \frac{r_{i'}}{2} \frac{\wp '(x)}{\wp (x) -\wp (\delta _{i'})}.
\end{equation}
By expanding Eq.(\ref{eq:llp}) at $x=0$ and observing coefficient of $x^0$, we obtain
\begin{equation}
c -\sum_{j=1}^l \zeta(t_j )= \left. \frac{\sqrt{-Q}}{\Xi (x)}\right| _{x=0},
\end{equation}
because the functions $\sigma '(x)/\sigma (x)$, $\sigma ' _i (x)/\sigma _i (x)$, $\wp '(x)/(\wp (x) -\wp (\delta _{i'}))$ and $\Xi '(x)/\Xi (x)$ are odd and $\sigma '(-t)/\sigma (-t)= -\zeta (t)$.  
It follows from $Q \neq 0$ and Proposition \ref{prop:char} that, if $l_0 \neq 0$, then $\left. \frac{\sqrt{-Q}}{\Xi (x)}\right| _{x=0} =0$, and if $l_0 = 0$, then $\left. \frac{\sqrt{-Q}}{\Xi (x)}\right| _{x=0} $ is finite. Thus we obtain (v).

We show (vi). 
The function $\Lambda (x) \Lambda (-x) $ is even doubly-periodic and satisfies Eq.(\ref{prodDE}), because it is a product of the solutions to Eq.(\ref{eq:H}). Since the dimension of the space of the solutions to Eq.(\ref{eq:H}), which are even doubly-periodic, is one, we have $\Xi (x)=C \Lambda (x) \Lambda (-x) $ for some non-zero constant $C$. Hence we have $\Xi (t_j ) = \Xi (-t_j )=0$. 
On the other hand, we have $\Lambda (-t_j) \neq 0$ from (ii).
At $x=-t_j$, the l.h.s. of Eq.(\ref{eq:llp}) is finite, and the denominator of the r.h.s. is zero. Therefore we have 
\begin{equation}
\Xi ' (x)| _{x=-t_j} +2\sqrt{-Q} =0.
\end{equation}
By changing the variable $z=\wp (x)$ and the oddness of the function $\wp '(x)$, we obtain (vi).
\end{proof}

\section{The case $M=1$, $r_1 =1$ and Painlev\'e equation} \label{sec:P6}
We consider Eq.(\ref{eq:Hg}) for the case $M=1$, $r_1 =1$.
For this case, Eq.(\ref{eq:Hg}) is written as 
\begin{equation}
(H_g-\tilde{E})f_g (x)=0,
\label{Hgkr1}
\end{equation}
where
\begin{align}
&  H_g= -\frac{d^2}{dx^2} + \frac{\wp ' (x)}{\wp (x) -\wp (\delta _{1})} \frac{d}{dx} + \frac{\tilde{s}_{1}}{\wp (x) -\wp (\delta _{1})} +\sum_{i=0}^3 l_i(l_i+1) \wp (x+\omega_i) . \label{HgP6}
\end{align}
We set
\begin{align}
& \Psi _g (x)= \sqrt{\wp (x) -\wp (\delta _1)} , \quad b_1 =\wp (\delta _1) , \label{pgkr0} \\
& \mu _1= \frac{-\tilde{s}_{1}}{4 b_1^3 -g_2 b_1 -g_3} +\sum _{i=1}^3 \frac{l_i}{2(b_1-e_i)}, \\
& p= \tilde{E} -2(l_1l_2e_3+l_2l_3e_1+l_3l_1e_2) +\sum _{i=1}^3 l_i (l_i e_i +2(e_i+b_1 )) .  \label{pgkr}
\end{align}
The condition that, the regular singular points $x=\pm \delta _1$ is apparent, is written as
\begin{align}
& p=(4 b_1^3 -g_2 b_1 -g_3) \left\{ -\mu _{1} ^2 +\sum_{i=1}^3\frac{l_i+\frac{1}{2}}{b_1-e_i} \mu _{1} \right\} \label{pgkrap} \\
& \quad \quad -b_1 (l_1+l_2+l_3-l_0)(l_1+l_2+l_3+l_0+1) .\nonumber
\end{align}
From now on we assume that $l_0, l_1, l_2, l_3 \in \Zint _{\geq 0}$ and the eigenvalue $\tilde{E}$ satisfies Eqs.(\ref{pgkr}, \ref{pgkrap}). Then the assumption in Proposition \ref{prop:prod} is true, and propositions and theorem in the previous section are valid. The function $\Xi (x)$ in Proposition \ref{prop:prod} is written as 
\begin{equation}
\Xi (x)=c_0+\frac{d_0}{(\wp (x)-\wp (\delta _1))}+\sum_{i=0}^3 \sum_{j=0}^{l_i-1} b^{(i)}_j \wp (x+\omega_i)^{l_i-j} .
\label{Fxkr1}
\end{equation}
It follows from Proposition \ref{prop:Xionedim} that the function $\Xi (x)$ is determined uniquely up to multiplicative constant.
Ratios of the coefficients $c_0/d_0$ and $b^{(i)}_j/d_0$ $(i=0,1,2,3, \: j=0,\dots ,l_i-1)$ are written as rational functions in variables $b_1$ and $\mu _{1} $, because the coefficients $b^{(i)}_j$, $c_0$ and $d_0$ satisfy linear equations whose coefficients are rational functions in $b_1$ and $\mu _1$, which are obtained by substituting Eq.(\ref{Fxkr1}) into Eq.(\ref{prodDE}). The value $Q$ is calculated by Eq.(\ref{const}) and it is expressed as a rational function in $b_1$ and $\mu _{1} $ multiplied by $d_0 ^2$. We set 
\begin{equation}
\Lambda ( x)=\sqrt{\Xi (x)}\exp \int \frac{ \sqrt{-Q}dx}{\Xi (x)}, \quad \Lambda _g( x) = \Lambda ( x) \Psi _g(x).
\label{integ1P6}
\end{equation}
Due to Proposition \ref{prop:Linteg}, the function $\Lambda _g (x)$ is a solution to the differential equation (\ref{Hgkr1}). 
By Theorem \ref{thm:alpha}, the eigenfunction $\Lambda _g (x)$ is also expressed in the form of the Hermite-Krichever Ansatz. Namely, it is expressed as 
\begin{align}
& \Lambda _g (x) = \exp \left( \kappa x \right) \left( \sum _{i=0}^3 \sum_{j=0}^{\tilde{l}_i-1} \tilde{b} ^{(i)}_j \left( \frac{d}{dx} \right) ^{j} \Phi _i(x, \alpha ) \right)
\label{LalphaP6}
\end{align}
or
\begin{align}
& \Lambda _g (x) = \exp \left( \bar{\kappa } x \right) \left( \bar{c} +\sum _{i=0}^3 \sum_{j=0}^{\tilde{l}_i-2} \bar{b} ^{(i)}_j \left( \frac{d}{dx} \right) ^{j} \wp (x+\omega _i) +\sum_{i=1}^3 \bar{c}_i \frac{\wp '(x)}{\wp (x)-e_i} \right)
\label{Lalpha0P6}
\end{align}
where $l= l_0+l_1+l_2+l_3 +1$, $\tilde{l} _0 = l_0 +1$ and $\tilde{l}_i =l_i$ $(i=1,2,3)$. Now we investigate the values $\alpha $ and $\kappa $ in Eq.(\ref{LalphaP6}). Note that, if $\alpha \not \equiv 0$ (mod $2\omega _1 \Zint \oplus 2\omega _3 \Zint $), then the function $\Lambda _g (x)$ is expressed as Eq.(\ref{LalphaP6}) and we have
\begin{align}
& \Lambda _g (x+2\omega _j) = \exp (-2\eta _j \alpha +2\omega _j \zeta (\alpha ) +2 \kappa \omega _j ) \Lambda _g (x) , \quad  (j=1,3). \label{ellint00} 
\end{align}
\begin{prop} \label{prop:P6HK}
Assume that $M=1$, $r_1 =1$, $l_0, l_1, l_2, l_3 \in \Zint _{\geq 0}$ and the value $p$ satisfies Eq.(\ref{pgkrap}). Let $\alpha $ and $\kappa $ be the values determined by the Hermite-Krichever Ansatz  (see Eq.(\ref{LalphaP6})).
Then $\wp (\alpha )$ is expressed as a rational function in variables $b_1$ and $\mu _{1} $, $\wp ' (\alpha )$ is expressed as a product of $\sqrt {-Q}$ and a rational function in variables $b_1$ and $\mu _{1} $, and $\kappa $ is expressed as a product of $\sqrt {-Q}$ and a rational function in variables $b_1$ and $\mu _{1} $.
\end{prop}
\begin{proof}
We assume that $Q \neq 0$. For the case $Q=0$, the proposition is shown by considering a continuation from the case $Q \neq 0$.

It follows from Eqs.(\ref{eq:tilL}, \ref{periods}, \ref{rel:sigmai}) that
\begin{align} 
& \Lambda _g (x+2\omega _j ) = \exp \left( 2\eta _j \left( - \sum_{i'=1}^l t_{i'} + \sum_{i=1}^3 l_{i}\omega_{i} \right)  +2 \omega _j \left( c - \sum _{i=1}^3 l_i \eta_i \right) \right) \Lambda _g (x)  \label{ellinttj}  
\end{align}
for $j=1,3$. By comparing with Eq.(\ref{ellint}), we have
\begin{align}
& -2\eta _1 \alpha +2\omega _1 (\zeta (\alpha ) + \kappa )= -2\eta _1 \left( \sum_{i'=1}^l t_{i'} - \sum_{i=1}^3 l_{i}\omega_{i}\right)  +2\omega _1 \left( c - \sum _{i=1}^3 l_i \eta_i \right) +2\pi \sqrt{-1}n_1 \label{eq:n1} ,\\
& -2\eta _3 \alpha +2\omega _3 (\zeta (\alpha ) + \kappa )= -2\eta _3 \left(\sum_{i'=1}^l t_{i'} - \sum_{i=1}^3 l_{i}\omega_{i}\right)  +2\omega _3 \left( c - \sum _{i=1}^3 l_i \eta_i \right)  +2\pi \sqrt{-1}n_3, \label{eq:n3}
\end{align}
for integers $n_1$, $n_3$. It follows that
\begin{align}
& \left( \alpha - \left(\sum_{i'=1}^l t_{i'} - \sum_{i=1}^3 l_{i}\omega_{i}\right) \right)(-2\eta_1 \omega _3 + 2\eta_3 \omega _1 ) = 2\pi \sqrt {-1} (n_1\omega _3 -n_3 \omega _1 ) , \label{eq:aln} \\
& \left( \zeta (\alpha ) + \kappa -c+ \sum _{i=1}^3 l_i \eta_i \right)(2\eta_3 \omega _1 -2\eta_1 \omega _3 ) = 2\pi \sqrt {-1} (n_1\eta _3 -n_3 \eta _1 ) . \label{eq:kan} 
\end{align}
From Legendre's relation $\eta_1 \omega _3 - \eta_3 \omega _1 = \pi \sqrt {-1}/2$, we have 
\begin{equation}
\alpha \equiv \sum_{i'=1}^l t_{i'} - \sum_{i=1}^3 l_{i}\omega_{i} \quad \quad  (\mbox{mod }2\omega_1 \Zint \oplus 2\omega_3 \Zint).
\end{equation}
Combining Eqs.(\ref{eq:aln}, \ref{eq:kan}) with Proposition \ref{prop:BA} (v) and relations $\zeta (\alpha +2\omega _i )= \zeta (\alpha ) +2\eta _i$ $(i=1,3)$, we have 
\begin{equation}
\kappa = -\zeta \left(\sum_{j=1}^{l} t_j - \sum_{i=1}^3 l_i \omega_i \right) + \sum_{j=1}^{l} \zeta (t_j)- \sum_{i=1}^3 l_i \eta_i + \delta _{l_0,0} \frac{\sqrt{-Q}}{\Xi (0)}.
\end{equation}

Next, we investigate values $\wp (\alpha )$, $\wp '(\alpha )$ and $\kappa $.
The functions $\wp (\sum_{j=1}^l t_j - \sum_{i=1}^3 l_{i}\omega_{i} )$, $\wp '(\sum_{j=1}^l t_j - \sum_{i=1}^3 l_{i}\omega_{i} )$ and 
$\zeta (\sum_{j=1}^{l} t_j - \sum_{i=1}^3 l_{i}\omega_{i}) - \sum_{j=1}^{l} \zeta (t_j) + \sum_{i=1}^3 l_{i}\eta_{i}$ are doubly-periodic in variables $t_1 ,\dots ,t_l$.
Hence by applying addition formulae of elliptic functions and considering the parity of functions $\wp (x)$, $\wp '(x)$ and $\zeta (x)$, we obtain the expression
\begin{align}
& \wp \left(\sum_{j=1}^l t_j - \sum_{i=1}^3 l_{i}\omega_{i}\right) = \sum_{j_1<j_2<\dots <j_m\atop{m:\mbox{\scriptsize{ even}}}}f^{(1)}_{j_1, \dots ,j_m}(\wp(t_1), \dots ,\wp (t_l)) \wp'(t_{j_1}) \dots \wp'(t_{j_l}), \\
& \wp '\left(\sum_{j=1}^l t_j - \sum_{i=1}^3 l_{i}\omega_{i}\right) = \sum_{j_1<j_2<\dots <j_m\atop{m:\mbox{\scriptsize{ odd}}}}f^{(2)}_{j_1, \dots ,j_m}(\wp(t_1), \dots ,\wp (t_l)) \wp'(t_{j_1}) \dots \wp'(t_{j_l}), \nonumber \\
& \zeta \left(\sum_{j=1}^{l} t_j - \sum_{i=1}^3 l_{i}\omega_{i}\right) - \sum_{j=1}^{l} \zeta (t_j) + \sum_{i=1}^3 l_{i}\eta_{i} \nonumber \\
& = \sum_{j_1<j_2<\dots <j_m\atop{m:\mbox{\scriptsize{ odd}}}}f^{(3)}_{j_1, \dots ,j_m}(\wp(t_1), \dots ,\wp (t_l)) \wp'(t_{j_1}) \dots \wp'(t_{j_l}), \nonumber
\end{align}
where $f^{(i)}_{j_1, \dots ,j_m}(x_1, \dots , x_l)$ $(i=1,2,3)$ are rational functions in $x_1 ,\dots ,x_l$.
From Eq.(\ref{signpptj}), the function  $\wp '(t_j)/\sqrt{-Q}$ is expressed as a rational function in $b_1$, $\mu _{1} $ and $\wp (t_j)$.
Hence, $\wp (\sum_{j=1}^l t_j - \sum_{i=1}^3 l_{i}\omega_{i})$, $\wp '(\sum_{j=1}^l t_j - \sum_{i=1}^3 l_{i}\omega_{i})/\sqrt{-Q}$ and $(\zeta (\sum_{j=1}^{l} t_j- \sum_{i=1}^3 l_{i}\omega_{i}) - \sum_{j=1}^{l} \zeta (t_j) + \sum_{i=1}^3 l_{i}\eta_{i})/\sqrt{-Q}$ are expressed as rational functions in the variable $\wp(t_1), \dots ,\wp (t_l)$, $b_1$ and $\mu _{1} $, and they are symmetric in $\wp(t_1), \dots ,\wp (t_l)$.

Since the dimension of the space of the solutions to Eq.(\ref{eq:H}), which are even doubly-periodic, is one, we have $\Xi (x)=C \Lambda (x) \Lambda (-x) $ for some non-zero scalar $C$.
Hence, we have the following expression;
\begin{equation}
\Xi (x) \Psi _g(x) ^2 =\frac{D \prod_{j=1}^{l}(\wp(x)-\wp(t_j))}{(\wp(x)-e_1)^{l_1}(\wp(x)-e_2)^{l_2}(\wp(x)-e_3)^{l_3}} \label{Fxtj}
\end{equation}
for some value $D (\neq 0)$. Thus
\begin{equation}
\prod_{j=1}^{l}(\wp(x)-\wp(t_j))= \Xi (x) \Psi _g (x) ^2 (\wp(x)-e_1)^{l_1}(\wp(x)-e_2)^{l_2}(\wp(x)-e_3)^{l_3}/D.
\end{equation}
Hence, the elementary symmetric functions $ \sum _{j_1<\dots <j_{l'}} \wp(t_{j_1})\dots \wp (t_{j_{l'}})$ ($l'=1,\dots ,l$) are expressed as rational functions in $b_1$ and $\mu _{1} $.
By substituting elementary symmetric functions into the symmetric expressions of $\wp (\sum_{j=1}^l t_j - \sum_{i=1}^3 l_{i}\omega_{i})$, $\wp (\sum_{j=1}^l t_j - \sum_{i=1}^3 l_{i}\omega_{i})$ and $(\zeta (\sum_{j=1}^{l} t_j- \sum_{i=1}^3 l_{i}\omega_{i}) - \sum_{j=1}^{l} \zeta (t_j) + \sum_{i=1}^3 l_{i}\eta_{i})/\sqrt{-Q}$, it follows that $\wp (\sum_{j=1}^l t_j - \sum_{i=1}^3 l_{i}\omega_{i})$, $\wp '(\sum_{j=1}^l t_j - \sum_{i=1}^3 l_{i}\omega_{i})/\sqrt{-Q}$ and $(\zeta (\sum_{j=1}^{l} t_j- \sum_{i=1}^3 l_{i}\omega_{i}) - \sum_{j=1}^{l} \zeta (t_j) + \sum_{i=1}^3 l_{i}\eta_{i})/\sqrt{-Q}$ are expressed as rational functions in $b_1$ and $\mu _{1} $.
Hence, $\wp (\alpha )$, $\wp '(\alpha )/ \sqrt{-Q}$ and $\kappa / \sqrt{-Q}$ are expressed as rational functions in variables $b_1$ and $\mu _{1} $.
\end{proof}

We now discuss the relationship between the monodromy preserving deformation of Fuchsian equations and the sixth Painlev\'e equation. For this purpose we recall some definitions and results of Painlev\'e equation.

The sixth Painlev\'e equation is a non-linear ordinary differential equation written as
\begin{align}
\frac{d^2\lambda }{dt^2} = & \frac{1}{2} \left( \frac{1}{\lambda }+\frac{1}{\lambda -1}+\frac{1}{\lambda -t} \right) \left( \frac{d\lambda }{dt} \right) ^2 -\left( \frac {1}{t} +\frac {1}{t-1} +\frac {1}{\lambda -t} \right)\frac{d\lambda }{dt} \label{eq:P6eqn} \\
& +\frac{\lambda (\lambda -1)(\lambda -t)}{t^2(t-1)^2}\left\{ \frac{\kappa _{\infty}^2}{2} -\frac{\kappa _{0}^2}{2}\frac{t}{\lambda ^2} +\frac{\kappa _{1}^2}{2}\frac{(t-1)}{(\lambda -1)^2} +\frac{(1-\kappa _{t}^2)}{2}\frac{t(t-1)}{(\lambda -t)^2} \right\}. \nonumber
\end{align}
A remarkable property of this differential equation is that its solutions do not have movable singularities other than poles.
This equation is also written in terms of a Hamiltonian system by adding the variable $\mu$, which is called the sixth Painlev\'e system:
\begin{equation}
\frac{d\lambda }{dt} =\frac{\partial H_{VI}}{\partial \mu}, \quad \quad
\frac{d\mu }{dt} =-\frac{\partial H_{VI}}{\partial \lambda}
\label{eq:Psys}
\end{equation}
with the Hamiltonian 
\begin{align}
H_{VI} = & \frac{1}{t(t-1)} \left\{ \lambda (\lambda -1) (\lambda -t) \mu^2 \right. \label{eq:P6} \\
& \left. -\left\{ \kappa _0 (\lambda -1) (\lambda -t)+\kappa _1 \lambda (\lambda -t) +(\kappa _t -1) \lambda (\lambda -1) \right\} \mu +\kappa (\lambda -t)\right\} ,\nonumber
\end{align}
where $\kappa = ((\kappa _0 +\kappa _1 +\kappa _t -1)^2- \kappa _{\infty} ^2)/4$.
The sixth Painlev\'e equation for $\lambda $ is obtained by eliminating $\mu $ in Eq.(\ref{eq:Psys}).
Set $\omega _1=1/2$, $\omega _3=\tau /2$ and write
\begin{equation}
t= \frac{e_3- e_1}{e_2-e_1}, \quad \lambda = \frac{\wp (\delta )- e_1}{e_2-e_1}.
\end{equation}
Then the sixth Painlev\'e equation is equivalent to the following equation (see \cite{Man,Tks}):
\begin{equation}
\frac{d^2 \delta }{d \tau ^2} = -\frac{1}{4\pi ^2} \left\{ \frac{\kappa _{\infty}^2}{2} \wp ' \left(\delta  \right) + \frac{\kappa _{0}^2}{2} \wp ' \left(\delta +\frac{1}{2} \right) + \frac{\kappa _{1}^2}{2} \wp ' \left(\delta +\frac{\tau +1}{2} \right) +  \frac{\kappa _{t}^2}{2} \wp ' \left(\delta +\frac{\tau }{2}\right) \right\}, \label{eq:P6ellip}
\end{equation}
where $\wp ' (z ) = (\partial /\partial z ) \wp (z)$.

It is widely known that the sixth Painlev\'e equation is obtained by the monodnomy preserving deformation of a certain linear differential equation.
Let us introduce the following Fuchsian differential equation:
\begin{equation}
\frac{d^2y}{dw^2} + p_1 (w) \frac{dy}{dw} +p_2(w) y=0, \label{eq:mpdP6}
\end{equation}
where 
\begin{align}
& p_1 (w) = \frac{1-\kappa _0}{w} + \frac{1-\kappa _1}{w-1} + \frac{1-\kappa _t}{w-t} -\frac{1}{w-\lambda}, \\
& p_2 (w) = \frac{\kappa }{w(w-1)} -\frac{t(t-1) H_{VI}}{w(w-1)(w-t)} + \frac{\lambda (\lambda -1) \mu}{w(w-1)(w-\lambda)}.
\end{align}
This equation has five regular singular points $\{ 0,1,t,\infty ,\lambda \}$ and the exponents at $w=\lambda $ are $0$ and $2$.
It follows from Eq.(\ref{eq:P6}) that the regular singular point $w=\lambda $ is apparent.
Then the sixth Painlev\'e equation is obtained by the monodromy preserving deformation of Eq.(\ref{eq:Psys}), i.e., the condition that the monodromy of Eq.(\ref{eq:mpdP6}) is preserved as deforming the variable $t$ is equivalent to that $\mu $ and $\lambda $ satisfy the Painlev\'e system (see Eq.(\ref{eq:Psys})), provided $\kappa _0, \kappa _1, \kappa _t , \kappa _{\infty} \not \in \Zint$. For details, see \cite{IKSY}.

Now we transform Eq.(\ref{eq:mpdP6}) into the form of Eq.(\ref{HgP6}). We set 
\begin{align}
& w=\frac{\wp (x) -e_1}{e_2-e_1}, \quad y= f_g (x) \prod _{i=1}^3 (\wp (x)-e_i)^{l_i/2}, \label{eq:wwpx} \\
&  \quad t=\frac{e_3-e_1}{e_2-e_1}, \quad \lambda =\frac{b_1 -e_1}{e_2 -e_1}, \quad \wp (\delta _1) =b_1.
\end{align}
Then we obtain Eq.(\ref{HgP6}) by setting
\begin{align}
& \kappa _0 =l_1 +1/2, \quad \kappa _1 =l_2 +1/2, \quad \kappa _t =l_3 +1/2, \quad \kappa _{\infty} =l_0 +1/2, \label{eq:kili} \\
& \mu = (e_2-e_1)\mu _1, \quad \kappa = (l_1+l_2+l_3+l_0 +1)(l_1+l_2+l_3-l_0) , \\
& H_{VI}=\frac{1}{t(1-t)} \left\{ \frac{p+\kappa e_3}{e_2-e_1} +\lambda (1- \lambda)\mu \right\},
\end{align}
(see Eqs.(\ref{pgkr0}--\ref{pgkr})), and Eq.(\ref{eq:P6}) is equivalent to Eq.(\ref{pgkrap}), that means that the apparency of regular singularity is inheritted.
Mapping from the variable $x$ to the variable $w$ (see Eq.(\ref{eq:wwpx})) is a double covering from the punctured torus $(\Cplx / (2\omega _1 \Zint + 2\omega _3 \Zint)) \setminus \{ 0, \omega_1 , \omega _2, \omega _3 \} $ to the punctured Riemann sphere ${\mathbb P}^1 \setminus \{ 0,1,t,\infty \}$. A solution $y(w)$ to Eq.(\ref{eq:mpdP6}) corresponds to a solution $f _g (x) $ to Eq.(\ref{HgP6}) by $y(w)= f_g (x) \prod _{i=1}^3 (\wp (x)-e_i)^{l_i/2}$. Hence the monodromy preserving deformation of Eq.(\ref{eq:mpdP6}) in $t$ corresponds to the monodromy preserving deformation of Eq.(\ref{HgP6}) in $\tau $.

Now we consider monodromy preserving deformation in the variable $\tau$ ($\omega _1 =1/2, \omega _3=\tau /2$) by applying solutions obtained by the Hermite-Krichever Ansatz for the case $l_i \in \Zint _{\geq 0}$ $(i=0,1,2,3)$.
Let $\alpha $ and $\kappa $ be values determined by the Hermite-Krichever Ansats (see Eq.(\ref{LalphaP6})). We consider the case $Q\neq 0$. Then a basis for solutions to Eq.(\ref{eq:Hg}) is given by $\Lambda _g(x)$ and $\Lambda _g(-x)$, and the monodromy matrix with respect to the cycle $x \rightarrow x+2\omega _j $ ($j=1,3$) is diagonal. The elements of the matrix are obtained from Eq.(\ref{ellint00}).
Hence, the eigenvalues $\exp (\pm ( -2\eta _j \alpha +2\omega _j \zeta (\alpha ) +2 \kappa \omega _j))$ $(j=1,3)$ of the monodromy matrices are preserved by the monodromy preserving deformation.
We set 
\begin{align}
& -2\eta _1 \alpha +2\omega _1 \zeta (\alpha ) +2 \kappa \omega _1 = \pi \sqrt{-1} C_1, \\
& -2\eta _3 \alpha +2\omega _3 \zeta (\alpha ) +2 \kappa \omega _3 = \pi \sqrt{-1} C_3, 
\end{align}
for contants $C_1$ and $C_3$.
By Legendre's relation, we have 
\begin{align}
& \alpha  = C_3 \omega _1 -C_1 \omega _3  \label{al00},\\
& \kappa = \zeta (C_1 \omega _3 -C_3 \omega _1 ) +C_3 \eta _1 -C_1 \eta _3  , \label{kapp00}
\end{align}
(see Eqs.(\ref{al}, \ref{kapp})). From Proposition \ref{prop:P6HK}, the value $\wp (\alpha)(=\wp (C_3 \omega _1-C_1 \omega _3 ))$ is expressed as a rational function in variables $b_1$ and $\mu _1$, the value $\wp '(\alpha)(=\wp '(C_3 \omega _1-C_1 \omega _3 ))$ is expressed as a product of $\sqrt {-Q}$ and a rational function in variables $b_1$ and $\mu _1$, and the value $\kappa (= \zeta (C_1 \omega _3 -C_3 \omega _1 ) +C_3 \eta _1 -C_1 \eta _3)$ is expressed as a product of $\sqrt {-Q}$ and rational function in variables $b_1$ and $\mu _1$. By solving these equations for $b_1$ and $\mu _1$ and evaluating them into Eq.(\ref{HgP6}), the monodromy of the solutions on the cycles $x \rightarrow x +2\omega _i$ $(i=1,3)$ are preserved for the fixed values $C_1$ and $C_3$. 
Let $\gamma _0$ be the path in the $x$-plane which is obtained by the pullback of the cycle turning the origin around anti-clockwise in the $w$-plane, where $x$ and $w$ are related with $w=(\wp (x) -e_1)/(e_2 -e_1)$.
Then the monodromy matrix on $\gamma _0$ with respect to the basis $(\Lambda _g(x) ,\Lambda _g(-x) )$ is written as 
\begin{equation}
(\Lambda _g(x) ,\Lambda _g(-x) ) \rightarrow (\Lambda _g(-x) ,\Lambda _g(x) ) = (\Lambda _g(x) ,\Lambda _g(-x) ) 
\left(
\begin{array}{cc}
0 & 1 \\
1 & 0
\end{array}
\right) ,
\end{equation}
and does not depend on $\tau $.
Since the fundamental group on the punctured Riemann sphere ${\mathbb P}^1 \setminus \{ 0,1,t,\infty \}$ is generated by the images of $\gamma _0$ and the cycles $x \rightarrow x +2\omega _i$ $(i=1,3)$, Eqs.(\ref{al00}, \ref{kapp00}) describe the condition for the monodromy preserving deformation on the punctured Riemann sphere by rewriting the variable $\tau$ to $t$.
Summarizing, we have the following proposition.
\begin{prop} \label{prop:P6}
We set $\omega _1= 1/2$, $\omega _3 =\tau /2$ and assume that $l_i \in \Zint _{\geq 0}$ $(i=0,1,2,3)$ and $Q\neq 0$.
By solving the equations in Proposition \ref{prop:P6HK} in variable $b_1 =\wp (\delta _1)$ and $\mu _1$, we express $\wp (\delta _1)$ and $\mu _1$ in terms of $\wp (\alpha)$, $\wp '(\alpha)$ and $\kappa $, and we replace $\wp (\alpha)$, $\wp '(\alpha)$ and $\kappa $ with $\wp (C_3 \omega _1-C_1 \omega _3 )$, $\wp '(C_3 \omega _1-C_1 \omega _3 )$ and $\zeta (C_1 \omega _3 -C_3 \omega _1) +C_3 \eta _1 -C_1 \eta _3$. Then $\delta _1$ satisfies the sixth Painlev\'e equation in the elliptic form
\begin{equation}
\frac{d^2 \delta _1}{d \tau ^2} = -\frac{1}{8\pi ^2} \left\{ \sum _{i=0}^3 (l_i +1/2)^2 \wp '( \delta _1 + \omega _i) \right\}. \label{eq:P6ellipl}
\end{equation}
\end{prop}
We observe the expressions of $b_1$ and $\mu _1$ in detail for the cases $l_0=l_1=l_2=l_3=0$ and $l_0=1$, $l_1=l_2=l_3=0$.

\subsection{The case $M=1$, $r_1 =1$, $l_0=l_1=l_2=l_3=0$}

We investigate the case $M=1$, $r_1 =1$, $l_0=l_1=l_2=l_3=0$ in detail.
The differential equation (\ref{Hgkr1}) is written as 
\begin{equation}
\left\{ -\frac{d^2}{dx^2} + \frac{\wp ' (x)}{\wp (x) -b_1} \frac{d}{dx} - \frac{ \mu _{1} (4 b_1^3 -g_2 b_1 -g_3)}{\wp (x) -b_1} -p \right\} f_g (x)=0,
\label{Hgkr1l00}
\end{equation}
We assume that $b_1 \neq e_1, e_2, e_3$. The condition that the regular singular points $x= \pm \delta _1$ $(\wp (\delta _{1})=b_1 )$ are apparent is written as
\begin{align}
& p=- (4 b_1^3 -g_2 b_1 -g_3) \mu _{1} ^2 +(6b_1^2 -g_2/2) \mu _{1}  \label{pgkrapl00} 
\end{align}
(see Eq.(\ref{pgkrap})). The doubly-periodic function $\Xi (x)$ (see Eq.(\ref{Fxkr1})) which satisfies Eq.(\ref{prodDE}) is calculated as 
\begin{equation}
\Xi (x)= 2\mu _1 +\frac{1}{\wp(x)-b_1} .
\end{equation}
The value $Q$ (see Eq.(\ref{const})) is calculated as
\begin{align}
& Q= 2\mu _1(2\mu _1 (e_1-b_1)+1)(2(e_2-b_1) \mu _1+1)(2\mu _1(e_3-b_1)+1). 
\end{align}
We set
\begin{equation}
\Lambda _g(x) = \sqrt{\Xi (x) (\wp (x) - b_1)} \exp \int \frac{ \sqrt{-Q}dx}{\Xi (x)},
\label{integ1P6l00}
\end{equation}
(see Eq.(\ref{integ1P6})).
Then a solution to Eq.(\ref{Hgkr1l00}) is written as $\Lambda _g (x)$, and is expressed in the form of the Hermite-Krichever Ansatz as
\begin{align}
& \Lambda _g (x) = \exp (\kappa x) \Phi _0 (x, \alpha )
\end{align}
for generic $(\mu_1 , b_1)$.
The values $\alpha $ and $\kappa $ are determined as
\begin{align}
& \wp (\alpha )= b_1 - \frac{1}{2 \mu_1}, \quad \wp '(\alpha )= -\frac{\sqrt{-Q}}{2\mu _1^2} , \quad \kappa = \frac{\sqrt{-Q}}{2\mu _1}.
\end{align}
Hence we have
\begin{align}
& \mu _1 = -\frac{\kappa  }{\wp ' (\alpha )} ,\quad  b_1 = \wp (\alpha ) -\frac{\wp ' (\alpha )}{2\kappa }.
\end{align}
From Proposition \ref{prop:P6}, the function $\delta _1$ determined by
\begin{align}
\wp (\delta _1) =  b_1 & = \wp (C_3 \omega _1-C_1 \omega _3 ) -\frac{\wp ' (C_3 \omega _1-C_1 \omega _3 )}{2(\zeta (C_1 \omega _3 -C_3 \omega_1 ) -C_1 \eta _3 +C_3 \eta _1) } \label{P6sol0000} \\
& = \wp (C_1 \omega _3 -C_3 \omega_1 ) +\frac{\wp ' (C_1 \omega _3 -C_3 \omega_1 )}{2(\zeta (C_1 \omega _3 -C_3 \omega_1 ) -(C_1 \eta _3 -C_3 \eta _1)) } \nonumber
\end{align}
is a solution to the sixth Painlev\'e equation in the elliptic form (see Eq.(\ref{eq:P6ellipl})). This solution coincides with the one found by Hitchin \cite{Hit} when he studied Einstein metrics and isomonodromy deformations.

Now we consider the case $Q=0$. If $Q=0$, then $\mu _1=0$ or $\mu _1 =1/(2(b_1 -e_i))$ for some $i \in \{1,2,3 \}$. 

If $\mu _1 =0$, then a solution to Eq.(\ref{Hgkr1l00}) is $1 (= \Lambda _g (x))$ and another solution is written as 
\begin{equation}
\textstyle \zeta (x) +b_1 x (=\int -(\wp (x) -b_1) dx ) .
\end{equation} 
We investigate the monodromy preserving deformation on the basis $s _1 (x)= B(\tau )$ and $s_2(x )=  \zeta (x) +b_1 x$, where $B(\tau )$ is a constant that is independent of $x$.
The monodromy matrix with respect to the path $\gamma _0$ is written as diag$(1,-1)$.
Since $s_2 (x + 2\omega _i ) = s_2 (x) +2(\eta _i + \omega _i b_1)$ $(i=1,3)$, the monodromy matrix with respect to the basis $(s_1 (x), s _2 (x))$ on the cycle $x \rightarrow x+2\omega _i$ $(i=1,3)$ is written as 
\begin{equation}
\left( 
\begin{array}{cc}
1 & \frac{2(\eta _i + \omega _i b_1)}{B(\tau )} \\
0 & 1
\end{array}
\right).
\end{equation}
To preserve monodromy, the matrix elements should be constants of the variable $\tau (=\omega _3/\omega _1)$. Hence we obtain 
\begin{align}
& 2(\eta _1 + \omega _1 b_1) =D_1 B(\tau ) ,\\
& 2(\eta _3 + \omega _3 b_1) =D_3 B(\tau ) , \nonumber
\end{align}
for some constants $D_1$ and $D_3$. By using Legendre's relation, we obtain that $ B(\tau ) = \pi \sqrt{-1} /(D_1 \omega _3 -D_1 \omega _1)$ and 
\begin{equation}
\wp (\delta _1) =b_1 =-\frac{D_1 \eta _3 -D_3 \eta _1}{D_1 \omega _3 -D_3 \omega _1}.
\label{b1mu10}
\end{equation}
Since Eq.(\ref{b1mu10}) is obtained by monodromy preserving deformation, the function $\delta _1$ satisfies the sixth Painlev\'e equation.

If $\mu _1 =1/(2(b_1 -e_i))$  for some $i \in \{1,2,3 \}$, then $\wp _i (x) (= \Lambda _g (x))$ is a solution to Eq.(\ref{Hgkr1l00}), and another solution is written as
\begin{equation}
\textstyle \wp _i (x) \left\{ \frac{e_i- b_1}{(e_i-e_{i'})(e_i-e_{i''})} \zeta (x +\omega _i)  +\right. \left. (1- \frac{e_i- b_1}{(e_i-e_{i'})(e_i-e_{i''})}) x \right\} (= \wp _i (x) \int \frac{\wp (x) -b}{\wp (x)-e_i} dx ),
\end{equation}
where $i'$ and $i''$ are elements in $\{ 1,2,3 \}$ such that $i' \neq i$, $i'' \neq i$ and $i' <i''$.
By calculating similarly to the case $\mu _1=0$, we obtain that the function $\delta _1$, which is determined by
\begin{equation}
\wp (\delta _1) =b_1 =\frac{(g_2/4 -2e_i^2)(D_1 \omega _3 -D_3 \omega _1) +e_i (D_1 \eta _3 -D_3 \eta _1) }{e_i (D_1 \omega _3 -D_3 \omega _1) +(D_1 \eta _3 -D_3 \eta _1) },
\label{b1mui}
\end{equation}
is a solution to the sixth Painlev\'e equation for constants $D_1$ and $D_3$.

We now show that Eqs.(\ref{b1mu10} ,\ref{b1mui}) are obtained by suitable limits from Eq.(\ref{P6sol0000}). Set $(C_1 ,C_3 )= (CD_1, CD_3)$ in Eq.(\ref{P6sol0000}) and consider the limit $C \rightarrow 0$, then we recover Eq.(\ref{b1mu10}). Similarly, set $(C_1 ,C_3) =(C D_1 ,-1+CD_3 )$ (resp. $(C_1 ,C_3) =(-1 +C D_1 ,1+ CD_3 )$, $(C_1 ,C_3) =(1+C D_1 , CD_3 )$) and consider the limit $C \rightarrow 0$, then we recover Eq.(\ref{b1mui}) for the case $i=1$ (resp. $i=2$, $i=3$).
Hence the space of the parameters of the solutions to the sixth Painlev\'e equation (i.e. the space of initial conditions) for the case $l_0=l_1=l_2=l_3=0$ is obtained by blowing up four points on the surface $\Cplx /(2\pi \sqrt{-1} \Zint ) \times \Cplx /(2\pi \sqrt{-1} \Zint ) $, and this reflects the $A_1 \times A_1 \times A_1 \times A_1$ structure of Riccati solutions by Saito and Terajima \cite{ST}.

\subsection{The case $M=1$, $r_1 =1$, $l_0=1$, $l_1=l_2=l_3=0$}
The differential equation (\ref{Hgkr1}) for this case is written as 
\begin{equation}
\left\{ -\frac{d^2}{dx^2} + \frac{\wp ' (x)}{\wp (x) -b_1} \frac{d}{dx} - \frac{ \mu _{1} (4 b_1^3 -g_2 b_1 -g_3)}{\wp (x) -b_1} + 2\wp (x) -p \right\} f_g (x)=0,
\label{Hgkr1l01}
\end{equation}
We assume that $b_1 \neq e_1, e_2, e_3$. The condition that the regular singular points $x= \pm \delta _1$ $(\wp (\delta _{1})=b_1 )$ are apparent is written as
\begin{align}
& p=- (4 b_1^3 -g_2 b_1 -g_3) \mu _{1} ^2 +(6b_1^2 -g_2/2) \mu _{1} +2b_1 \label{pgkrapl01} 
\end{align}
(see Eq.(\ref{pgkrap})). The doubly-periodic function $\Xi (x)$ (see Eq.(\ref{Fxkr1})), which satisfies Eq.(\ref{prodDE}), is calculated as 
\begin{align}
\Xi (x)=& \wp (x) +( (-4b_1^3+b_1g_2+g_3)\mu _1^2+(6b_1^2-g_2/2)\mu _1 -b_1) \\
& +  ((-4b_1^3+b_1g_2+g_3)\mu _1 /2+3b_1^2-g_2/4)/(\wp (x) -b_1) \nonumber
\end{align}
The value $Q$ (see Eq.(\ref{const})) is calculated as
\begin{align}
 Q= -& ((2(4 b_1^3-b_1g_2-g_3)\mu _1^3-(12b_1^2-g_2)\mu _1^2+4) (2(b_1^2+e_1b_1+e_2e_3)\mu _1-2b_1-e_1) \\
& (2(b_1^2+e_2b_1+e_1e_2)\mu _1-2b_1-e_2) (2(b_1^2+e_3b_1+e_1e_3)\mu _1-2b_1-e_3) \nonumber .
\end{align}
We set
\begin{equation}
\Lambda _g(x) = \sqrt{\Xi (x) (\wp (x) - b_1)} \exp \int \frac{ \sqrt{-Q}dx}{\Xi (x)}, 
\label{integ1P6l01}
\end{equation}
(see Eq.(\ref{integ1P6})).
Then a solution to Eq.(\ref{Hgkr1l00}) is written as $\Lambda _g (x)$, and it is expressed in the form of the Hermite-Krichever Ansatz as
\begin{align}
& \Lambda _g (x) = \exp (\kappa x) \left\{ \Phi _0 (x, \alpha ) +\frac{d}{dx} \Phi _0 (x, \alpha ) \right\} 
\end{align}
for generic $(\mu_1 , b_1)$.
The values $\alpha $ and $\kappa $ are determined as
\begin{align}
& \wp (\alpha )  = \frac{2(4b_1^3-b_1g_2-g_3)b_1 \mu _1^3+(-24b_1^3+4g_2b_1+3g_3)\mu _1^2+(24b_1^2-2g_2)\mu _1-8b_1}{2(4b_1^3-b_1 g_2-g_3)\mu _1^3 -(12b_1^2-g_2)\mu _1^2+4},\\
& \wp '(\alpha )= \frac{-4((4b_1^3-b_1g_2-g_3)\mu _1^3-(12b_1^2-g_2)\mu _1^2+12b_1\mu _1-4)}{(2(4b_1^3-b_1 g_2-g_3)\mu _1^3 -(12b_1^2-g_2)\mu _1^2+4)^2}\sqrt{-Q}  ,\\
& \kappa = \frac{2\mu _1}{2(4b_1^3-b_1 g_2-g_3)\mu _1^3 -(12b_1^2-g_2)\mu _1^2+4}\sqrt{-Q}.
\end{align}
Hence we have
\begin{align}
& b_1 = \frac{2\wp (\alpha ) \kappa ^3-3\wp '(\alpha )\kappa ^2+(6\wp (\alpha ) ^2 -g_2)\kappa -\wp (\alpha ) \wp '(\alpha )}{2(\kappa ^3-3\wp (\alpha ) \kappa +\wp '(\alpha ))} ,\\
& \mu _1 = \frac{2(\kappa ^3-3\wp (\alpha ) \kappa +\wp '(\alpha ))\kappa }{-2\wp '(\alpha )\kappa ^3+(12\wp (\alpha ) ^2-g_2)\kappa ^2-6\wp (\alpha ) \wp '(\alpha )\kappa +\wp '(\alpha )^2}.
\end{align}
From Proposition \ref{prop:P6}, the function $\delta _1$ determined by
\begin{align}
& \wp (\delta _1) = b_1 = \label{P6sol1000} \\ 
& \quad  \frac{2\wp (\omega ) (\zeta (\omega )- \eta )^3+3\wp '(\omega )(\zeta (\omega )- \eta )^2+(6\wp (\omega ) ^2 -g_2)(\zeta (\omega )- \eta )+\wp (\omega ) \wp '(\omega )}{2((\zeta (\omega )- \eta )^3-3\wp (\omega ) (\zeta (\omega )- \eta ) -\wp '(\omega ))} , \nonumber \\
& (\omega = C_1 \omega _3 -C_3 \omega_1 , \quad  \eta = C_1 \eta _3 -C_3 \eta _1) ,\nonumber
\end{align}
is a solution to the sixth Painlev\'e equation in the elliptic form (see Eq.(\ref{eq:P6ellipl})).
In the sixth Painlev\'e equation, it is known that the case $(\kappa _{0}, \kappa _{1}, \kappa _{t}, \kappa _{\infty}) =(1/2, 1/2, 1/2, 3/2 ) $ is linked to the case $(\kappa _{0}, \kappa _{1}, \kappa _{t}, \kappa _{\infty}) =(1/2, 1/2, 1/2, 1/2 ) $ by B\"acklund transformation. 
For a table of B\"acklund transformation of the sixth Painlev\'e equation, see \cite{TOS}.
By transformating the solution in Eq.(\ref{P6sol0000}) of the case $(\kappa _{0}, \kappa _{1}, \kappa _{t}, \kappa _{\infty}) =(1/2, 1/2, 1/2, 1/2 )$ to the one of the case  $(\kappa _{0}, \kappa _{1}, \kappa _{t}, \kappa _{\infty}) =(1/2, 1/2, 1/2, 3/2 )$, we recover the solution in Eq.(\ref{P6sol1000}).

Now we consider the case $Q=0$. If $Q=0$, then $\mu_1  $ is a solution to the equation $2(4 b_1^3-b_1g_2-g_3)\mu _1^3-(12b_1^2-g_2)\mu _1^2+4 =0$ or $\mu _1=(2b_1+e_i)/(2(b_1^2+e_ib_1+e_i^2-g_2/4))$ for some $i \in \{1,2,3 \}$.
We set $\omega = D_1 \omega _3 -D_3 \omega_1$ and $\eta = D_1 \eta _3 -D_3 \eta _1$, where $D_1$ and $D_3$ are constants. For the case that $\mu_1  $ is a solution to the equation $2(4 b_1^3-b_1g_2-g_3)\mu _1^3-(12b_1^2-g_2)\mu _1^2+4 =0$, the corresponding solutions to the sixth Painlev\'e equation are written as the function $\delta _1$, where 
\begin{equation}
\wp( \delta _1)= b_1 = \frac{4 \eta ^3 +g_2 \omega ^2 \eta -2 g_3 \omega ^3}{\omega (g_2 \omega ^2-12\eta ^2 )}. \label{b1mu1l01}
\end{equation}
For the case $\mu _1=(2b_1+e_i)/(2(b_1^2+e_ib_1+e_i^2-g_2/4))$ ($i \in \{1,2,3 \}$), we have 
\begin{equation}
\wp( \delta _1)= b_1 = \frac{-g_2 e_i \omega /2 +(6e_i^2 -g_2)\eta}{(6e_i^2 -g_2)\omega -6e_i\eta }. \label{b1muil01}
\end{equation}
Note that these solutions are also obtained by suitable limits from Eq.(\ref{P6sol1000}), and Eq.(\ref{b1mu1l01}) (resp. Eq.(\ref{b1muil01})) is transformed by B\"acklund transformation from Eq.(\ref{b1mu10}) (resp. Eq.(\ref{b1mui})).

\section{Relationship with finite-gap potential} \label{sec:FG}

\subsection{The case $M=0$ and Heun's equation}
For the case $M=0$, Eq.(\ref{Feq}) is transformed to Heun's equation, and the potential of the operator $H(=-\frac{d^2}{dx^2} + v(x))$ (see Eq.(\ref{Ino})) is written as 
\begin{equation}
v(x)= \sum_{i=0}^3 l_i(l_i+1)\wp (x+\omega_i) .\label{HHeun}
\end{equation}
If $l_0, l_1, l_2, l_3 \in \Zint _{\geq 0}$, then the function $v(x)$ in Eq.(\ref{HHeun}) is called the Treibich-Verdier potential, and is an example of algebro-geometric finite-gap potential (see \cite{TV,GW,Smi,Tak3}). Recall that, if there exists an odd-order differential operator
\begin{equation}
A= \left( \frac{d}{dx} \right)^{2g+1} + \sum_{j=0}^{2g-1}  b_j(x) \left( \frac{d}{dx} \right)^{2g-1-j}
\end{equation}
 such that 
\begin{equation}
\left[A, -\frac{d^2}{dx^2}+v(x)\right]=0,
\end{equation}
 then $v(x)$ is called an algebro-geometric finite-gap potential.
For the case $M=0$ and $l_0, l_1, l_2, l_3 \in \Zint _{\geq 0}$, statements in Propositions \ref{prop:prod}, \ref{prop:Linteg}, \ref{prop:indep}, \ref{prop:char}, \ref{prop:zeros}, \ref{prop:BA} and Theorem \ref{thm:alpha} in this paper hold true, because there is no constraint relation for the apparency of additional regular singularity.
Moreover, the function $\Xi (x)$ in Proposition \ref{prop:prod} is written as 
\begin{equation}
\Xi (x)= E^g +\sum _{i=1}^{g} a_{g-i} (x) E^{i},
\label{XiHeun}
\end{equation}
for some $g \in \Zint _{\geq 1}$ and even doubly-periodic functions $a _{i} (x)$ $(i=0, \dots ,g-1)$, and the constant $Q$ is a monic polynomial in $E$ of degree $2g +1$ (see \cite{Tak1}).
The commuting operator $A$ is described by using functions $a _{i}(x)$ $(i=1, \dots ,g)$.
The eigenfunction $\Lambda (x)$ (see Eq.(\ref{integ1})) of the operator $-\frac{d^2}{dx^2}+v(x)$ is expressed in a form of the Hermite-Krichever Ansatz (see Theorem \ref{thm:alpha}). It is shown in \cite{Tak4} that the values $\wp (\alpha )$, $\wp ' (\alpha )/\sqrt{-Q}$ and $\kappa /\sqrt{-Q}$ are expressed as a rational function in $E$, and it follows that the global monodromy of Heun's equation for the case $l_0, l_1, l_2, l_3 \in \Zint _{\geq 0}$ is written as an elliptic integral.
On the other hand it is known that global monodromy is also expressed by a hyperelliptic integral (see \cite{Tak3}). By comparing the two expressions, we obtain a hyperelliptic-to-elliptic integral reduction formula (see \cite{Tak4}).

\subsection{The case $M=1$ and $r_1=2$}

For the case $M=1$ and $r_1=2$, the differential equation in the elliptic form is written as 
\begin{equation}
\left( -\frac{d^2}{dx^2} + v(x)-E\right) f(x)= 0,
\label{InoEFk1}
\end{equation}
where $v(x)$ is written as 
\begin{align}
v(x) = & 2(\wp (x-\delta _{1}) + \wp (x+\delta _{1})) +\frac{s_{1}}{\wp (x) -\wp (\delta _{1})} +\sum_{i=0}^3 l_i(l_i+1)\wp (x+\omega_i) . \label{Inok1}
\end{align}
Set $b_1 =\wp (\delta _1)$. The condition that the regular singularity $x =\pm \delta _{1} $ of Eq.(\ref{InoEFk1}) is apparent (which is equivalent to that the regular singularity $z= b_1$ of Eq.(\ref{Ino}) is apparent) is written as 
\begin{align}
& s_1^3+ (12b_1^2- g_2)s_1^2  +(4(4b_1^3-g_2 b_1 -g_3)E+ f_1(b_1))s_1  +f_0(b_1) =0.
\end{align}
where $f_1(b_1) $ and $f_0(b_1)$ are given by
\begin{align}
f_1(b_1) =& -2(2l_0^2+2l_0+5)b_1 (4b_1^3-g_2 b_1 -g_3) +(6b_1 ^2 -g_2/2)^2 \\
& -8(2l_1^2+2l_1+1)(b_1-e_2)(b_1-e_3)(e_1 b_1+e_1^2+e_2e_3 ) \nonumber \\
& -8(2l_2^2+2l_2+1)(b_1-e_1)(b_1-e_3)(e_2 b_1+e_2^2+e_1e_3 ) \nonumber \\
& -8(2l_3^2+2l_3+1)(b_1-e_1)(b_1-e_2)(e_3 b_1+e_3^2+e_1e_2 ) ,\nonumber \\
f_0(b_1) =& (2l_0+1)^2(4b_1^3-g_2 b_1 -g_3)^2 \\
& - 16(2l_1+1)^2 (e_1-e_2)(e_1-e_3)(b_1-e_2)^2(b_1-e_3)^2 \nonumber \\
& - 16(2l_2+1)^2 (e_2-e_1)(e_2-e_3)(b_1-e_1)^2(b_1-e_3)^2  \nonumber \\
& - 16(2l_3+1)^2 (e_3-e_1)(e_3-e_2)(b_1-e_1)^2(b_1-e_2)^2 . \nonumber 
\end{align}
If $s_1 =0$, then we obtain an equation
\begin{align}
& f_0(b_1) =0.
\label{eq:f0b10}
\end{align}
Then the value $b_1$ is determined by this equation, and does not depend on the value $E$.

Assume that $s_1 =0$ and the value $b_1$ satisfies  Eq.(\ref{eq:f0b10}). Then the function $\Xi (x)$ is expressed as Eq.(\ref{XiHeun}), and similar properties are valid as in similar arguments written on Heun's equation in \cite{Tak3}.
In fact, a commuting operator of an odd degree is constructed from the function $\Xi (x)$, and we recover the results by Treibich \cite{Tre}. Moreover, a solution to Eq.(\ref{InoEFk1}) is expressed in the form of the Hermite-Krichever Ansatz, monodmony has two integral representations that are ellptic and hyperelliptic, and hyperelliptic-to-elliptic integral reduction formulae are obtained. Details will be reported in \cite{TakF}.

Thus, if we restrict our discussion to the case $s_1 =0$ and that $b_1$ satisfies Eq.(\ref{eq:f0b10}), then the potential of the operator $H$ (see Eq.(\ref{Inok1})) is a finite-gap, which recover the results by Treibich \cite{Tre} and Smirnov \cite{Smi2}. In other words, the potential is Picard's in the sense of Gesztesy and Weikard \cite{GW2}. Note that Smirnov obtained expressions like Eqs.(\ref{XiHeun}, \ref{integ1}) and calculated several examples in \cite{Smi2}.

\section{Concluding remarks} \label{sec:rmk}

We have shown in sections \ref{sec:HK} and \ref{sec:P6} that solutions of the linear differential equation that produces the sixth Painlev\'e equation have integral representations and that they are expressed in the form of the Hermite-Krichever Ansatz. Furthermore we got a procedure for obtaining solutions of the sixth Painlev\'e equation (see Eq.(\ref{eq:P6ellip})) for the cases $\kappa _0 , \kappa _1 , \kappa _t, \kappa _{\infty} \in \Zint +\frac{1}{2}$ by fixing the monodromy, and we presented explicit solutions for the cases $(\kappa _0 , \kappa _1 , \kappa _t, \kappa _{\infty} )=(\frac{1}{2}, \frac{1}{2}, \frac{1}{2}, \frac{1}{2})$ and $(\frac{1}{2}, \frac{1}{2}, \frac{1}{2}, \frac{3}{2})$.

By B\"acklund transformation of the sixth Painlev\'e equation (see \cite{TOS} etc.), Hitchin's solution (i.e., solutions for the case $(\kappa _0 , \kappa _1 , \kappa _t, \kappa _{\infty} )=(\frac{1}{2}, \frac{1}{2}, \frac{1}{2}, \frac{1}{2})$) is transformed to the solutions for the case $(\kappa _0 , \kappa _1 , \kappa _t, \kappa _{\infty} ) \in O_1 \cup O_2$, where
\begin{align}
& O_1= \left\{ (\kappa _0 , \kappa _1 , \kappa _t, \kappa _{\infty} ) | 
\kappa _0 , \kappa _1 , \kappa _t, \kappa _{\infty} \in \Zint +\frac{1}{2} \right\}, \\
& O_2 = \left \{(\kappa _0 , \kappa _1 , \kappa _t, \kappa _{\infty} ) \left| 
\begin{array}{ll}
\kappa _0 , \kappa _1 , \kappa _t, \kappa _{\infty} \in \Zint \\
\kappa _0 + \kappa _1 + \kappa _t + \kappa _{\infty}  \in 2 \Zint 
\end{array}
\right. \right\}. 
\end{align}
Note that solutions for the case $(\kappa _0 , \kappa _1 , \kappa _t, \kappa _{\infty} )=(0, 0, 0, 0)( \in O_2)$ are already known and are called Picard's solution.

For the case $(\kappa _0 , \kappa _1 , \kappa _t, \kappa _{\infty} ) \in O_1$, solutions of the linear differential equation are investigated by our method, and solutions of the sixth Painlev\'e equation follow from them.
On the other hand, for the case $(\kappa _0 , \kappa _1 , \kappa _t, \kappa _{\infty} ) \in O_2$, we cannot obtain results on integral representation and the Hermite-Krichever Ansatz by our method, although solutions of the sixth Painlev\'e equation are obtained in principle by B\"acklund transformation.
Note that the condition $(\kappa _0 , \kappa _1 , \kappa _t, \kappa _{\infty} ) \in O_1$ corresponds to the condition $l_0, \dots ,l_3 \in \Zint +\frac{1}{2}$, $l_0+ l_1 +l_2 + l_3 \in 2 \Zint $.

Now we propose a problem to investigate solutions and their monodromy of the linear differential equation (Eq.(\ref{Hgkr1}) with the condition (\ref{pgkrap})) for the cases $l_0, \dots ,l_3 \in \Zint +\frac{1}{2}$, $l_0+ l_1 +l_2 + l_3 \in 2 \Zint $. In partiuclar, how can we investigate solutions and their monodromy of the linear differential equation for the case $\kappa _0 = \kappa _1 = \kappa _t= \kappa _{\infty} =0$ (i.e. $l_0=l_1=l_2=l_3 =-1/2$)?

\appendix
\section {Elliptic functions} \label{sect:append}
This appendix presents the definitions of and the formulas for the elliptic functions.

The Weierstrass $\wp$-function, the Weierstrass sigma-function and the Weierstrass zeta-function with periods $(2\omega_1, 2\omega_3)$ are defined as follows:
\begin{align}
& \wp (z)= \frac{1}{z^2}+  \sum_{(m,n)\in \Zint \times \Zint \setminus \{ (0,0)\} } \left( \frac{1}{(z-2m\omega_1 -2n\omega_3)^2}-\frac{1}{(2m\omega_1 +2n\omega_3)^2}\right),  \\
& \sigma (z)=z\prod_{(m,n)\in \Zint \times \Zint \setminus \{(0,0)\} } \left(1-\frac{z}{2m\omega_1 +2n\omega_3}\right) \nonumber \\
& \; \; \; \; \; \; \; \; \; \; \; \; \; \; \cdot \exp\left(\frac{z}{2m\omega_1 +2n\omega_3}+\frac{z^2}{2(2m\omega_1 +2n\omega_3)^2}\right), \nonumber \\
& \zeta(z)=\frac{\sigma'(z)}{\sigma (z)}. \nonumber
\end{align}
Setting $\omega_2=-\omega_1-\omega_3$ and 
\begin{align}
& e_i=\wp(\omega_i), \; \; \; \eta_i=\zeta(\omega_i) \; \; \; \; (i=1,2,3)
\end{align}
yields the relations
\begin{align}
& e_1+e_2+e_3=\eta_1+\eta_2+\eta_3=0, \; \; \; \label{eq:Leg} \\
& \eta _1 \omega _3- \eta _3 \omega _1 = \eta _3 \omega _2- \eta _2 \omega _3 = \eta _2 \omega _1- \eta _1 \omega _2 = \pi\sqrt{-1} /2, \nonumber \\
& \wp(z)=-\zeta'(z), \; \; \; (\wp'(z))^2=4(\wp(z)-e_1)(\wp(z)-e_2)(\wp(z)-e_3), \nonumber \\
& \wp (z) - \wp (\tilde{z}) = -\frac{\sigma ( z+\tilde{z})\sigma ( z-\tilde{z})}{\sigma ( z)^2\sigma (\tilde{z})^2} \nonumber
\end{align}
The periodicity of functions $\wp(z)$, $\zeta (z)$ and $\sigma (z)$ are as follows:
\begin{align}
& \wp(z+2\omega_i)=\wp(z), \; \; \; \zeta(z+2\omega_i)=\zeta(z)+2\eta_i \; \; \; \; (i=1,2,3), \label{periods} \\
& \sigma (z+2\omega _i) = - \sigma (z) \exp (2\eta _i (z + \omega _i)), \; \; \; \frac{\sigma (z+t+2\omega _i )}{\sigma (z+2\omega _i)}= \exp(2\eta _i t) \frac{\sigma (z+t)}{\sigma (z)} \nonumber
\end{align}
The constants $g_2$ and $g_3$ are defined by
\begin{equation}
g_2=-4(e_1e_2+e_2e_3+e_3e_1), \; \; \; g_3=4e_1e_2e_3.
\end{equation}
The co-sigma functions $\sigma_i(z)$ $(i=1,2,3)$ and co-$\wp$ functions $\wp_i(z)$ $(i=1,2,3)$ are defined by
\begin{align}
& \sigma_i(z)=\exp (-\eta_i z)\frac{\sigma(z+\omega_i)}{\sigma(\omega _i)}, \; \; \; \wp_i(z) = \frac{\sigma_i(z)}{\sigma(z)}, \label{eq:sigmai}
\end{align}
and satisfy
\begin{align}
& \wp_i(z) ^2 =\wp(z)-e_i, \quad \quad \quad (i,i' =1,2,3) \label{rel:sigmai} \\
& \wp _i (z+2\omega _{i'}) = \exp (2(\eta _{i'} \omega _i -\eta _{i} \omega _{i'}) ) \wp _i (z) = (-1)^{\delta _{i,i'}} \wp _i (z). \nonumber
\end{align}

{\bf Acknowledgments.}
The author would like to thank Professor Hidetaka Sakai for fruitful discussions. 
He is partially supported by a Grant-in-Aid for Scientific Research (No. 15740108) from the Japan Society for the Promotion of Science.


\begin{thebibliography}{9999}
\bibitem{BE}
Belokolos, E. D. and Enolskii, V. Z., Reduction of Abelian functions and algebraically integrable systems. II. 
{\it J. Math. Sci. (New York)} {\bf 108} (2002), 295--374.
\bibitem{GW}
Gesztesy F. and Weikard R., Treibich-Verdier potentials and the stationary (m)KdV hierarchy. {\it Math. Z.} {\bf 219} (1995), 451--476. 
\bibitem{GW2}
Gesztesy F. and Weikard R., Picard potentials and Hill's equation on a torus. {\it Acta Math.} {\bf 176} (1996), no. 1, 73--107.
\bibitem{Hit}
Hitchin N. J., Twistor spaces, Einstein metrics and isomonodromic deformations. {\it J. Differential Geom.} {\bf 42} (1995), no. 1, 30--112.
\bibitem{Inc}
Ince E. L. Ordinary Differential Equations. Dover Publications, New York, 1944.
\bibitem{IKSY}
Iwasaki K., Kimura H., Shimomura S. and Yoshida M., From Gauss to Painleve. A modern theory of special functions. Aspects of Mathematics, E16. Friedr. Vieweg \& Sohn, Braunschweig, 1991.
\bibitem{Man}
Manin, Yu. I., Sixth Painleve equation, universal elliptic curve, and mirror of $\bold P\sp 2$. Geometry of differential equations, 131--151, Amer. Math. Soc. Transl. Ser. 2, 186, Amer. Math. Soc., Providence, RI, 1998.
\bibitem{Ron}
Ronveaux A.(ed.), {\it Heun's differential equations.} Oxford Science Publications, Oxford University Press, Oxford, 1995.
\bibitem{ST}
Saito M.-H. and Terajima H.: Nodal curves and Riccati solutions of Painlev\'e equations. Preprint, math.AG/0201225. 
\bibitem{Smi}
Smirnov A. O.,  Elliptic solitons and Heun's equation, 
{\it The Kowalevski property,} 287--305, CRM Proc. 
Lecture Notes, 32, Amer. Math. Soc., Providence (2002).
\bibitem{Smi2}
Smirnov A. O.,  Finite-gap solutions of the Fuchsian equation, Preprint, math.CA/0310465, 2003. 
\bibitem{Tks}
Takasaki K., Painleve-Calogero correspondence revisited. {\it J. Math. Phys.} {\bf 42} (2001), no. 3, 1443--1473.
\bibitem{Tak1}
Takemura K., The Heun equation and the Calogero-Moser-Sutherland system I: the Bethe Ansatz method. {\it Comm. Math. Phys.} {\bf 235} (2003), 467--494.
\bibitem{Tak2}
Takemura K., The Heun equation and the Calogero-Moser-Sutherland system II: the perturbation and the algebraic solution, {\it  Electron. J. Differential Equations} {\bf 2004} no. 15 (2004), 1--30. 
\bibitem{Tak3}
Takemura K., The Heun equation and the Calogero-Moser-Sutherland system III: the finite gap property and the monodromy, {\it J. Nonlinear Math. Phys.} {\bf 11} (2004) 21--46.
\bibitem{Tak4}
Takemura K., The Heun equation and the Calogero-Moser-Sutherland system IV: the Hermite-Krichever Ansatz, Preprint, math.CA/0406141, 2004.
\bibitem{TakF}
Takemura K., in preparation.
\bibitem{Tre}
Treibich A. Hyperelliptic tangential covers, and finite-gap potentials.
{\it Russian Math. Surveys} {\bf 56} (2001), no. 6, 1107--1151.
\bibitem{TV}
Treibich A. and Verdier J.-L., Revetements exceptionnels et sommes de 4 nombres triangulaires (French). {\it Duke Math. J.} {\bf 68} (1992), 217--236. 
\bibitem{TOS}
Tsuda T., Okamoto K. and Sakai H.: Folding transformations of the Painlev\'e equations, Preprint UTMS 2003-42, 2003.

\end{thebibliography}
\end{document}